

\input amstex
\documentstyle{amsppt}
\magnification\magstep1

\topmatter
\title Nonnormal del Pezzo surfaces\endtitle
\author Miles Reid\endauthor
\affil Univ\. of Warwick\endaffil
\endtopmatter

\TagsOnRight

\hyphenation{Theo-rem dual-ising pre-dual-ising quasi-projective}

\define\simby#1{\buildrel{{\scriptstyle\text{#1}}}\over\sim} 
\define\bir{\simby{{bir}}} 
\define\myno#1{\par\vskip-\parskip\medskip\noindent{\rm#1\enspace}} 
\define\dual{\mathrel{\raise3pt\hbox{$\underline{\roman{\thinspace d
\thinspace}}$}}} 
\define\red{_{\hbox{\eightrm red}}} 
\define\gen{_{\hbox{\eightrm gen}}} 
\define\hidot{^{\textstyle\cdot}} 
\define\QED{\ifhmode\unskip\nobreak\fi\quad{Q.E.D.}}
\define\al{\alpha}
\define\de{\delta}
\define\De{\Delta}
\define\ep{\varepsilon}
\define\fie{\varphi}
\define\w{\omega}
\define\Ga{\Gamma}
\define\aff{\Bbb A} 
\define\Z{\Bbb Z} 
\define\FF{\Bbb F} 
\define\proj{\Bbb P} 
\define\Cee{{\Cal C}} 
\define\F{{\Cal F}} 
\define\G{{\Cal G}} 
\define\I{{\Cal I}} 
\define\J{{\Cal J}} 
\define\Oh{{\Cal O}} 
\define\CDer{{\Cal{D}}er} 
\define\Ext{\operatorname{Ext}} 
\define\CExt{{\Cal{E}}xt} 
\define\Hom{\operatorname{Hom}} 
\define\CHom{{\Cal{H}}om} 
\define\Ann{\operatorname{Ann}} 
\define\Ass{\operatorname{Ass}} 
\define\Res{\operatorname{Res}} 
\define\dd{\roman{d}} 
\define\Tr{\operatorname{Tr}} 
\define\im{\operatorname{im}} 
\define\Pico{\operatorname{Pic^0}} 
\define\Sing{\operatorname{Sing}} 
\define\Spec{\operatorname{Spec}} 
\define\card{\operatorname{card}} 
\define\cha{\operatorname{char}} 
\define\rest#1{_{\textstyle{\vert}#1}} 
\define\wave{\widetilde}

\heading\S0. Introduction\endheading

\myno{0.1} Throughout this paper, a {\it del Pezzo surface} is by definition a
connected, $2$@-dimensional, projective $k$@-scheme $X,\Oh_X(1)$ that is
Gorenstein and anti\-canonically polar\-ised; in other words, $X$ is
Cohen--Macaulay, and the dualising sheaf is invertible and anti\-ample: $\w_X
\cong\Oh_X(-1)$. For example, $X=X_3\subset\proj^3$ an arbitrary hypersurface
of degree $3$.

Under extra conditions, del Pezzo surfaces are interesting for lots of reasons:
for example, as tangent cones to index $1$ canonical $3$@-fold singularities
\cite{C3-f, 2.13}; by the ``Serre correspondence'', they could occur as the
subschemes in $\proj^4$ corresponding to sections of unstable vector bundles
over $\proj^4$ (the nonexistence proof of \cite{Grauert and Schneider} has a
gap). The main motivation for the present study was \cite{Mori,
Proposition~3.9}, where the statement that an irreducible del Pezzo surface $X$
has $\chi(\Oh_X)\ne 0$ plays an essential role; I reprove here, in particular,
Mori's statement that an irreducible del Pezzo surface in characteristic $0$
has $\chi(\Oh_X)=1$ (see Corollary~4.10). This was originally proved by
S.~Mori and S.~Goto (unpublished), and Mori was also kind enough to correct an
imbecility in my proof.

\myno{0.2} I assume throughout that $X$ is reduced, but either reducible or
(irreducible and) nonnormal; $\pi\:Y\to X$ is the normalisation. $Y$ has $r\ge
1$ components, and is marked by an effective Weil divisor $C$,
scheme-theoretically defined by the conductor ideal
$I_{C,Y}=\Cee_Y\subset\Oh_Y $; as a set, $C$ is the codimension $1$ double
locus of $\pi$. By {\it subadjunction} (Proposition~2.3), it's well known that
the canonical Weil divisor of $Y$ is $K_Y=\pi^*K_X-C$, so that
 $$
-K_Y=\pi^*\Oh_X(1)+C=(\text{ample})+(\text{effective}).
$$
It's easy to classify the components $C\subset Y$ with this property
(Theorem~1.1) in terms of scrolls and $\proj^2$, and in particular each
connected component of $C,\Oh_C(1)$ must be isomorphic to a plane conic.

Write $D\subset X$ for the subscheme defined by the conductor ideal
 $$
I_{D,X}=\Cee=\Ann(\pi_*\Oh_Y/\Oh_X)\subset\Oh_X,
 $$
and $\fie\:C\to D$ for the restriction of $\pi$. Thus $X$ is obtained by
glueing together one or more components $C_i\subset Y_i$ along a morphism
$\fie\:C=\coprod C_i\to D$. The Cohen--Macaulay or $S_2$ condition for $X$ is
easy: $\fie_*\Oh_C/\Oh_D$ must have no sections supported at points
Proposition~2.2, so that the glueing is entirely determined in codimension
$1$. The combinatorics of the glueing also turns out to be straightforward
(see 1.3 and Lemma~4.1); in particular if $Y$ is reducible then all the conics
$C_i$ are isomorphic.

\myno{0.3} What makes a local ring Gorenstein? We've not heard the last of this
question. The technical crux of this paper is Theorem~2.6, which characterises
the predualising sheaf $\w_X$ of a nonnormal scheme $X$ as the sheaf of {\it
Rosenlicht differentials}, that is, rational sections of the predualising
sheaf $\w_Y$ of the normalisation $Y$ with poles along $C$, and whose residues
along $C$ have zero trace on the generic fibres of $\fie\:C\to D$. This gives
the necessary and sufficient conditions of Corollary~2.8 on $C\subset Y$ and
$\fie\:C\to D$ for $X$ to be Gorenstein, for example: $\w_Y(C)$ is an
invertible $\Oh_Y$@-module, $\w_X$ is Gorenstein in codimension $1$ and $\w_C$
has an $\Oh_C$@-basis $s\in\ker\{\Tr_{C/D}\:\fie_*\w_C\to\w_D\}$. In
codimension $1$, where $\pi_*\Oh_Y$ local\-ised at a prime divisor of $X$ is a
product of DVRs, the condition for $\Oh_X\subset\pi_*\Oh_Y$ to be Gorenstein
translates into an interesting question on subrings $\Oh_D$ ``half-filling'' a
product $\fie_*\Oh_C=\prod A_i$ of Artin\-ian quotients of DVRs (see
3.4--5). This result may be useful in other contexts, since Gorenstein in
codimension $1$ is the essential prerequisite for working with conditions such
as log canonical singularities on nonnormal schemes. Even in the classic case
of curves over an algebraically closed field (Rosenlicht, Serre), it gives
rise to lots of unsolved problems, and my elementary results Theorem~3.7 may
be new even in this case.

The material on duality and Rosenlicht differentials in \S\S2--3 is written up
in much more generality than needed for the geometry of del Pezzo surfaces; it
is part of a continuing attempt to write up Grothendieck duality in an
absolutely elementary way (see also \cite{Reid}). \S3 is a laundry job on
\cite{Serre, Ch\.~IV, \S11}. (I hope that's not {\it l\`ese-majest\'e}.) I
haven't tried to find historically correct attributions for these ideas.

\myno{0.4} Working with the Gorenstein condition is easy if all the $C_i$ are
reduced, when $X$ has ordinary double points in codimension $1$ (see
Theorem~3.7, (I) and 4.2--3). But in the nonreduced case, the question is
quite subtle, and the answer depends on $\cha k$. In any characteristic,
$\Sing X=\Ga$ is an irreducible curve of degree $1$ for the polarisation, and
by Theorem~3.7, (II) the transverse singularity of $X$ along the generic point
of $\Ga$ is a cusp $(y^2=x^3)$ if $r=1$, a tac\-node $(y^2=x^2y)$ if $r=2$,
and $r$ concurrent lines in $\aff^{(r-1)}$ with no $(r-1)$ in a hyperplane for
$r\ge 3$ (if $r=3$, a plane ordinary triple point $xy(x-y)=0$). If $\cha k =
0$ then $\Ga\cong\proj^1$ and $D,\Oh_D(1)$ is isomorphic to a first order
neighbourhood of $\proj^1$ in $\proj^r$, with the morphism $\fie\:\coprod
C_i\to D$ linear on each component (see Proposition~4.9); this is essentially
equivalent to the main result $\chi(\Oh_X)=1$ (see Corollary~4.10, (I)).
However, in characteristic $p$ the curve $\Ga$ can have ``wild'' cusps (see
4.4, 4.7 and 4.11--12), essentially because the unknown subring
$\Oh_D\subset\fie_*\Oh_C$ is specified by a derivation (Proposition~3.9), and
$\roman dx^p\equiv0$. Any number of wild cusps can occur, and $\chi(\Oh_X)$
can be arbitrarily negative. Thus for nonnormal varieties in characteristic
$p$, Gorenstein is a weaker condition than in characteristic $0$, a new kind of
pathology discovered by Mori and Goto.

\myno{0.5} It follows easily that if $\cha k=0$ then $H^1(X,\Oh_X(n))=0$ for
all $n$, in particular $\chi(\Oh_X)=1$, and a general element $x_0\in
H^0(X,\Oh_X(1))$ is a non-zerodivisor for $\Oh_X$ (see Corollary~4.10, (I));
well-known arguments then show (Corollary~4.10, (II)) that $X$ has the usual
projective embedding properties of nonsingular del Pezzo surfaces or curves of
genus $1$: $\Oh_X(1)$ is very ample if $\deg X\ge 3$, the image is
ideal-theoretically an intersection of quadrics if $\deg X\ge 4$, and $X$ is a
weighted hypersurface or complete intersection $X_6\subset\proj(1^2,2,3)$,
$X_4\subset\proj(1^3,2)$, $X_3\subset \proj^3$ or $X_{2,2}\subset\proj^4$ if
$\deg X\le 4$. See 1.3--4 for examples.

However, if $\cha k=p$ then $H^1(X,\Oh_X)$ can be arbitrarily large; for an
irreducible surface $X$ this happens only when $X$ has a curve $\Ga$ of cusps,
and $\Ga$ itself has ``wild'' cusps. In particular, this cannot happen in
characteristic $\ge5$ if $X$ has only hypersurface singularities (see 4.4 and
4.7); I guess that a similar conclusion holds under other reasonable
conditions, e.g., $X$ lifts to characteristic $p^2$.

\proclaim{{\rm0.6}\enspace Conjecture} Assume that $\cha k=0$, and that $X$
is $1$@-connected but not reduced; then $X,\Oh_X(1)$ is projectively
Cohen--Macaulay (that is, $H^i(X,\Oh_X(n))=0$ for $i=1$ and all $n$, and
for $i=2$ and $n\ge0$), and a general element $x_0\in H^0(X,\Oh_X(1))$ is
$\Oh_X$@-regular, so all the projective embedding properties hold as above.
{\rm Finding the correct notion of $1$@-connected in general is part of the
problem; for example, if $X$ is locally a divisor in a smooth $3$@-fold and
$X=A+B$ then $\Oh_A(-B)$ is an invertible sheaf, and $1$@-connected means that
$\deg_{\Oh_X(1)}\Oh_A(-B)<0$. By adjunction, $\w_A=\w_X\otimes\Oh_A(-B)$, so
that, as in the reduced case, this amounts to $K_A{}^{-1}>(\text{ample
Cartier})$.}\endproclaim

\subheading{{\rm0.7}\enspace The great debate} Rend\. del circolo matematico di
Palermo {\bf1} (1887), p.~382 records the admission to the circle of {\it
dottore Pasquale del Pezzo, marchese di Campo\-disola}. It would be interesting
to know why Corrado Segre writing in the same volume (p.~218, 220, 221), along
with every subsequent Italian writer, spells the Marquis' name incorrectly
with a capital D.

\subheading{{\rm0.7}\enspace Acknowledgements} I thank S.~Mori for his help,
and for encouraging me to write up my proof, and M.~Miyanishi whose 1981
seminars on \cite{Mori} introduced me to the problem (and to much else
besides). The initial version of this paper was intended for the proceedings
of the wonderful Kinosaki conference in Dec~1981, and I apologise for the
delay in writing up the talk; however, as with other products consumed at
Kinosaki, the material undoubtedly improves on maturing for several years. I'm
very grateful to S.~Mori and K.~Saito for inviting me to Japan, and to RIMS,
Kyoto Univ\. for employing me during 1989--1990, when this paper was written.
J.~Koll\'ar and N.~Shepherd-Barron helped correct an error in the adjunction
formula of Theorem~2.12. Exercise~4.12 answers a question raised by
Shepherd-Barron and Koll\'ar.

\heading Contents\endheading

\S1. The normalised variety $C\subset Y$

\S2. Normalisation and dualising sheaves

\quad Appendix to \S2. Adjunction for a finite morphism

\S3. What happens in codimension $1$

\S4. The glueing map $\fie\:C\to D$ and proof of Theorem~1.5

References

\goodbreak

\heading\S1. The normalised variety $C\subset Y$\endheading

In (1.1--2), $Y$ denotes one component of the normalisation of $X$. As
explained in Theorem~2.3, it follows from subadjunction that $-K_Y=H+C$, where
$H=\Oh_Y(1)$ is an ample Cartier divisor and $C>0$ an effective Weil divisor. I
tabulate the information on $C\subset Y$ in the following form, to enable the
many readers familiar with the result to skip to the next section.

\proclaim{{\rm1.1}\enspace Theorem} Pairs $C\subset Y$ are listed as
follows:{\rm
 $$
\matrix\format\c\quad&\l\qquad&\c\qquad&\l\quad&\l\\
\text{Case}&Y,\Oh_Y(1)&\kern-1.25cm\text{degree $(\Oh_Y(1))^2$}\kern-.4cm
&\text{Class of $C$}&\text{Nature of $C$}\\
\vspace{12pt}
\text{(a)}&\proj^2,\Oh_{\proj^2}(1)&1&
\matrix\format\l\\
 \Oh_{\proj^2}(2),\text{that is,}\\ \vspace{3pt}
 \text{plane conic}\\
\endmatrix
&
\matrix\format\l\\
 \text{(a1) smooth conic} \\ \vspace{1pt}
 \text{(a2) line pair} \\ \vspace{1pt}
 \text{(a3) double line} \\
\endmatrix\\ \vspace{12pt}
\text{(b)}&\proj^2,\Oh_{\proj^2}(2)&4&\Oh_{\proj^2}(1)&\text{smooth conic}\\
\vspace{12pt}
\text{(c)}&
\matrix\format\l\\
 \FF_{a;0},aA\\ \vspace{3pt}
 \qquad\text{for $a\ge2$}\\
\endmatrix
&a&
\matrix\format\l\\
 2A,\text{ that is,}\\ \vspace{1pt}
 \text{2 generators,}\\ \vspace{1pt}
 \quad\text{or if $a=2$ only:}\\
\endmatrix
&
\matrix\format\l\\
 \text{(c1) line pair}\\ \vspace{1pt}
 \text{(c2) double line}\\ \vspace{1pt}
 \text{(c0) smooth conic}\\
\endmatrix\\ \vspace{12pt}
\text{(d)}&
\matrix\format\l\\
 \FF_{a;1},(a+1)A+B\\ \vspace{3pt}
 \qquad\text{for $a\ge0$}\\
\endmatrix
&a+2&
\matrix\format\l\\
 A+B\\ \vspace{1pt}
 \quad\text{or if $a\le1$ only:}\\
\endmatrix
&
\matrix\format\l\\
 \text{(d1) line pair}\\ \vspace{1pt}
 \text{(d0) smooth conic}\\
\endmatrix\\ \vspace{12pt}
\text{(e)}&\FF_{a;2},(a+2)A+B&a+4&B&\text{smooth conic}\\
\endmatrix
 $$
Here $\FF_a$ is the usual C.~Segre--P.~del Pezzo scroll \cite{Segre, del
Pezzo~1}, $A$ and $B$ its fibre and negative section; and
 $$
\FF_{a;k}=\bigl(\FF_a,\Oh((a+k)A+B)\bigr)
 $$
is the embedded scroll $\FF_{a;k}\subset\proj^{a+2k+1}$ of degree $a+2k$, with
negative section $B$ of degree $k$, except for $k=0$, when $\FF_{a;0}$ is the
cone over a rational normal curve, polarised by $\Oh(1)=\Oh(aA)$. The two
exceptional cases (c0) and (d0) with $a=0$ both correspond to a quadric of
$\proj^3$ with a smooth hyperplane section. Note that a double generator of a
cone is isomorphic to a double line in $\proj^2$, for example because it is a
projective cone over a tangent vector $k[\ep]/(\ep^2)$, or because
$\w_C\cong\Oh_C(-1)$.
 }\endproclaim

\demo{{\rm1.2} Proof} This theorem is presumably traditional, but Mori theory
gives an amazingly clean proof. I do a minimal resolution $f\:S\to Y$, set
$L=f^*H=f^*(-K_Y-C)$ and $C'=f'C$ (birational transform); then $L$ is a nef and
big
Cartier divisor, having positive intersection number with every component of
$C'$. Write
 $$
K_S+C'=f^*(K_Y+C)-Z;
 $$
clearly $Z$ is a Cartier divisor supported on the exceptional locus of $f$, and
$(K_S+C')\Ga\ge0$ for every exceptional curve $\Ga$, hence $Z\ge0$. Now
$K_S+L=-C'-Z$ is not nef, so that by the theorem on the cone \cite{Mori, 1.4
and 2.1} it follows that $S$ has an extremal rational curve $\ell$ with
$(K_S+L)\ell<0$. But $\ell$ cannot be a $(-1)$@-curve, since a $(-1)$@-curve
$m$ satisfies $K_Sm=-1$, $Lm>0$; therefore either $S=\proj^2$ and $\deg L\le
2$, or $S$ is a $\proj^1$@-bundle and $\deg L\le 1$ on the fibres. The theorem
follows
on sorting out the last possibility. \QED\enddemo

\subheading{{\rm1.3}\enspace Reassembling the pieces}
A nonnormal del Pezzo surface $X$ is obtained by glueing together a number
$r\ge1$ of the building blocks $C_i\subset Y_i$ of Theorem~1.1. How to ensure
that $X$ is Gorenstein is the main theme of \S\S2--4; since
$\Oh_Y(-K_Y-C)=\Oh_Y(1)$ is invertible, Corollary~2.8 reduces this to a
question on the glueing morphism $\fie\:C=\coprod C_i\to D$.

Here I present some classes of examples of the finished product $X$ in the
spirit of the projective geometry of the 1880s. 1.4 deals with cases specific
to $\deg X\le2$, when $\Oh_X(1)$ is not very ample.

\myno{(A)} Project $\FF_{a;2}\subset\proj^{a+5}$ to $\proj^{a+4}$ from a point
$P$ in the plane of the conic $B$ but not on $B$; the projection has a line
$\ell$ of ordinary double points. The same construction for the Veronese
surface
(case (b) of the theorem) is of course the well-known Bordello surface
$F^4_{(1)}[4]\subset\proj^4$ of \cite{Semple and Roth, p.~132}, which is a
complete intersection of two quadrics with a double line.

Note that in characteristic $2$, two cases occur, since all tangents to a plane
conic pass through a point $Q$; the projection $\fie\:B\to\ell$ from $Q$ is
inseparable, whereas the projection from any other point not on $B$ is
separable.

\myno{(B)} Let $C_i\subset Y_i$ be any two elements of (a1, b, c0, d0, e), not
both from (a1); that is, the $C_i$ are smooth conics, and not both of the $Y_i$
are planes. Embed $Y_1$ and $Y_2$ into a common projective space $\proj^n$ such
that the subspaces spanned by $Y_1$ and $Y_2$ intersect in the planes of the
$C_i$, and the $C_i$ are identified.

\myno{(C$_1$)} Let $C\subset Y$ be an element of (c1) or (d1) of degree
$\ge3$, so that $C$ is a line pair spanning a plane $\Pi$ of the projective
space $\proj^n$ containing $Y$ with $n\ge4$; make a linear projection $Y\to
X\subset\proj^{n-1}$ from a point $P\in\Pi\setminus C$. In case (c1), $X$ is
the projective cone over a nodal rational curve of degree $\ge3$. This case
can be viewed as the degenerate case $r=1$ of the following (C$_r$).

\myno{(C$_r$)} Let $C_i\subset Y_i$ be any $r\ge2$ elements of (a2, c1, d1),
not
both from (a2) if $r=2$; that is, the $C_i$ are line pairs, and not both
the $Y_i$ are planes if $r=2$. Embed the $Y_i$ into a common projective space
as a cycle of surfaces meeting along lines of $C_i$, with a common vertex.

In other words, name the lines of the pair $C_i$ as $C_i=\ell_i\cup\ell'_i$
(the
subscripts are cyclic, so $i=r+1$ counts as $i=1$). In $\proj^n$, choose a
vertex $P$ and $r$ linearly independent lines $m_i$ through $P$, and embed the
$Y_i$ so that $\ell_i$ is glued to $m_{i-1}$, and $\ell'_i$ to $m_i$, with the
subspaces spanned by $Y_i$ and $\coprod_{j\ne i} Y_j$ intersecting only in the
plane of $m_{i-1}$ and $m_i$, so that they are as linearly independent as
possible. The superfluous notation $\ell''_i=\ell_{i+1}$ will be convenient
later, so that the two lines of $C$ lying above $m_i$ are $\ell'_i$ and
$\ell''_i$.

\myno{(D$_1$)} The projective cone over a cuspidal rational curve of degree
$\ge3$. This can be viewed as a degenerate case of the following (D$_r$).

\myno{(D$_r$)} Let $C_i\subset Y_i$ be any $r\ge2$ elements of (a3, c2), not
both from (a3) if $r=2$; that is, the $C_i$ are double lines, and not both the
$Y_i$ are planes if $r=2$. Embed the $Y_i$ into a common projective space
meeting along a line, with a single nondegenerate linear dependence relation
between the $r$ planes of $C_i$.

In other words, in $\proj^n$, choose a line $\ell$ and $r$ planes $\pi_i$
through
$\ell$ such that any $(r-1)$ of them span the same $\proj^r$. Embed the $Y_i$
so
that $C_i$ is identified with the double line $2l\subset\pi_i$, and so that the
subspaces spanned by $Y_i$ and $\coprod_{j\ne i} Y_j$ intersect in $\pi_i$
only.

\subheading{{\rm1.4}\enspace Degree 1 and 2} I describe briefly the del Pezzo
surfaces of degree $1$ and $2$, leaving most of the computations to the
reader; a similar description of surfaces of degree $3$, $4$ and $5$ by
equations is also possible, and is an interesting exercise.

Consider first degree $1$. Let $Y=\proj^2$ with coordinates $u_1$, $u_2$,
$u_3$,
and $C:(q=0)$, where
 $$
q=u_2{}^2-u_1u_3,\quad u_2{}^2-u_1{}^2\quad\text{or}\quad u_2{}^2.
 $$
It's easy to see that setting $x_1=u_1$, $x_2=u_3$, $y=q$ and $z=u_2q$ defines
a birational morphism $\pi\:Y\to X_6\subset\proj(1^2,2,3)$, with image $X=X_6$
the sextic hypersurface defined by
 $$
z^2=y^3+x_1x_2y^2,\quad z^2=y^3+x_1{}^2y^2\quad\text{or}\quad z^2=y^3.
 $$
In characteristic $\ne2$, the map $\fie\:C\to D=\proj^1:(y=z=0)$ is the
quotient
by the $\Z/2$ action $u_2\mapsto -u_2$. In the second and third case, $X$ is a
weighted cone over a nodal or cuspidal rational curve, defined by
a polarisation of degree $1$.

In characteristic $2$, the same equations still define a normalisation, but in
the first case $\fie\:C\to D$ is an inseparable cover of $\proj^1$ by a
nonsingular conic; a separable example is provided by the normalisation of
$z^2+x_1yz=y^3+x_2{}^2y^2$.

It's interesting to observe that whereas for del Pezzo surfaces with isolated
singularities, those of degree $1$ are by far the most complicated, in the
nonnormal case those of degree $1$ are very few, and simple to describe.

Now for degree $2$. It's obvious that if $(q=q_2(x_1,x_2,x_3)=0)\subset\proj^2$
is a conic then the weighted quartic $X=X_4:(y^2=q^2)\subset\proj(1^3,2)$ is a
del Pezzo surface of degree $2$, consisting of $2$ copies of $Y\cong\proj^2$
glued along the conic $(q=0)$.

If $Y$ is irreducible of degree 2, then from the table of Theorem~1.1 it must
be
a quadric of $\proj^3$, with $C$ a hyperplane section. Consider, in
characteristic $\ne2$, the plane quartics $\ell^2q$ with a double line $\ell$
and $q$ a conic in the following $5$ cases: a nonsingular conic not tangent to
$\ell$, a nonsingular conic tangent to $\ell$, a line pair with vertex not on
$\ell$, a line pair with vertex on $\ell$ but not containing $\ell$, or a line
pair with $\ell$ as a component. Let $X=X_4:(y^2=\ell^2q)\subset\proj(1^3,2)$
be
the double cover of $\proj^2$ branched in $\ell^2q$. It's easy to see in the 5
cases that the normalisation of $X$ is a surface $C\subset Y$ of type (d0, d1,
c0, c1, c2).

\proclaim{{\rm1.5}\enspace Main theorem} Under the ``tame'' assumption of
4.7, the examples of 1.3--4 provide a complete classification of
nonnormal del Pezzo surfaces.\endproclaim

The ``tame'' condition covers all cases with $\cha k=0$, or $C$ reduced, or
$\cha k\ge5$ and $X$ locally a divisor in a nonsingular $3$@-fold; it is
equivalent to $H^1(X,\Oh_X)=0$ or $\chi(\Oh_X)=1$. See 4.4 and 4.11--12 for
counterexamples in the remaining cases.

\heading\S2. Normalisation and dualising sheaves\endheading

\myno{2.0} In this section $X$ is a purely $n$@-dimensional variety, and
$\pi\:Y\to X$ its normalisation. More generally $X$ could be a reduced
Noetherian scheme under the assumptions that the normalisation $\pi\:Y\to X$ is
finite and $Y$ has a dualising complex. The most general category for the
dualising complex has never been definitively established, but, for example,
if $Y$ is contained in a Gorenstein ambient scheme then dimensions and
codimensions are well defined (``universally catenary''), and $Y$ has a
dualising complex, whose top cohomology is the predualising sheaf $\w_Y$.

The aim is to use Grothendieck duality to describe $\w_Y$ in terms of $\w_X$
and vice versa. Note that all the sheaves here are coherent $\Oh_X$@-modules.

\myno{2.1} Write $\Cee=\Ann(\pi_*\Oh_Y/\Oh_X)\subset\Oh_X$ and
$\Cee_Y=\Cee\cdot\Oh_Y\subset\Oh_Y$ for the conductor of the normalisation, and
let $D\subset X$ and $C\subset Y$ be the subschemes they define. Then it
follows from the exact diagram
 $$
\spreadmatrixlines{3pt}
\matrix\format &\,\c\,\\
0&\to&\Cee&\subset&\Oh_X&\to&\Oh_D&\to&0\\
&&\Vert&&\bigcap&&\bigcap\\
0&\to&\pi_*\Cee_Y&\subset&\pi_*\Oh_Y&\to&\fie_*\Oh_C&\to&0\\
\endmatrix
\tag{$*$}
 $$
that
 $$
\pi_*\Oh_Y/\Oh_X=\fie_*\Oh_C/\Oh_D
\qquad\text{and}\qquad
\Cee=\Ann_{\Oh_X}(\fie_*\Oh_C/\Oh_D).
 $$
Note that by definition of $\Cee$ it follows that $\pi_*\Oh_Y/\Oh_X$ is a
faithful $\Oh_D$@-module.

All this just means that $X$ is the ringed space constructed by glueing a
normal
variety $Y$ along a finite morphism $\fie\:C\to D$; that is, $X$ is the
topological space $Y$ modulo the equivalence defined by $\fie$, and
 $$
\Oh_X=\ker\{\pi_*\Oh_Y\to\fie_*\Oh_C\to\fie_*\Oh_C/\Oh_D\}.
 $$
In general, starting from $Y$ and $\fie$, the resulting ringed space is not
necessarily a quasiprojective scheme (for example, the union of exceptional
surfaces in Hironaka's famous counterexample, or the fibre bundle of cuspidal
cubic curves of \cite{Horrocks}); is it always an algebraic space? In my case
this will never be a problem, since $\w_X{}^{-1}$ will be ample.

\proclaim{{\rm2.2}\enspace Proposition} $X$ satisfies $S_2$ at a
scheme-theoretic point $P\in X$ of codimension (height) $\ge2$ if and only if
every rational section of\/ $\Oh_X$ regular along every irreducible
codimension $1$ subscheme through $P$ is regular at\/ $P$; or, in other words,
for every scheme-theoretic point\/ $Q\in X$ of codimension $\ge2$ such that\/
$P\in V(Q)$ is in the closure of\/ $Q$, the $\Oh_X$@-submodule
$\Oh_X\subset\Oh_Y$ does not have a $Q$@-primary component, that is,
$Q\notin\Ass(\pi_*\Oh_Y/\Oh_X) =\Ass(\fie_*\Oh_C/\Oh_D)$.\endproclaim

\demo{Proof} The first sentence is standard, see for example \cite{YPG,
3.17--18}; since $Y$ is normal, a rational section of $\Oh_X$ regular along
every irreducible codimension $1$ subscheme through $P$ extends as a regular
section of $\pi_*\Oh_Y$, so that the $Q$@-primary part of $\pi_*\Oh_Y/\Oh_X$ is
exactly the obstruction to $\text{depth}_Q\,\Oh_X\ge2$.\enddemo

The proposition means the following: think of $X$ as constructed by glueing a
normal $n$@-dimensional variety $Y$ along a finite morphism $\fie\:C\to D$.
Then for $X$ to satisfy $S_2$, the glueing must all be forced by what happens
in codimension $1$; that is, $C$ and $D$ have pure codimension $1$, and
$\Oh_X\subset\pi_*\Oh_Y$ is determined by the local subrings
$\Oh_{D,\Ga}\subset(\fie_*\Oh_C)_\Ga$ at the generic point of each component
$\Ga$ of $D$, so that a given $f\in\pi_*\Oh_Y$ belongs to $\Oh_X$ if and only
if its image $\overline f\in\fie_*\Oh_C$ belongs to $\Oh_{D,\Ga}$ at the
generic point of each $\Ga$.

Note that if $C$ is reduced and $\fie\:C\to D$ separable, it follows that the
glueing is entirely geometric in nature, that is, $\Oh_X$ consists of all
functions on $Y$ constant on the general geometric fibres of $\fie$. See
4.2--3 for use of these ideas.

The $S_2$ condition (saturated in codimension $1$, reflexive, divisorial,
etc\.) is discussed for example in \cite{YPG, 3.17--18}. I make constant use
of the fact that if $\Oh_X$ is $S_2$ then $\Oh_X\subset\pi_*\Oh_Y$ is an
intersection of codimension $1$ primary components, and hence so is
$\Cee=\CHom(\pi_*\Oh_Y,\Oh_X)$. Thus each of $\Oh_X/\Cee$, $\pi_*\Oh_Y/\Cee$
and $\pi_*\Oh_Y/\Oh_X$ are $S_1$ or torsion-free as $\Oh_D$@-modules. This
condition can be interpreted as saying that $\Oh_D$ is normal under $\Oh_C$,
that is, rational sections of $\Oh_D$ in $\fie_*\Oh_C$ are already in $\Oh_D$.
Similar remarks apply to the duals $\w_D$, $\fie_*\w_C$ and $\ker\Tr_{C/D}$
(see Remark~2.9), so that many questions reduce to the stalks of these sheaves
at generic points of $D$.

The {\it $S_2$isation} or {\it saturation} of a coherent sheaf $\F$ on a
scheme $Z$ is the unique sheaf $\F'$ with an $\Oh_Z$@-linear morphism
$\F\to\F'$ such that $\F'$ is $S_2$. This is the same thing as the reflexive
hull or double dual of $\F$ if $Z$ is a normal variety.

 \proclaim\nofrills{{\rm2.3}\enspace Proposition}{\rm\ (``Subadjunction'',
compare \cite{Mumford}).\ }\ $\w_Y$ is determined by $\w_X$ and the conductor
as follows:
 $$
\pi_*\w_Y=\CHom_{\Oh_X}(\pi_*\Oh_Y,\w_X);
 $$
(the $\pi_*\Oh_Y$@-module structure is given by multiplication in the first
entry of the $\CHom$). The right-hand side is the biggest $\pi_*\Oh_Y$@-module
contained in $\w_X$. It coincides with the $S_2$isation of the $\Oh_Y$@-module
$\Cee\cdot\w_X$ if\/ $X$ is Gorenstein in codimension $1$ {\rm (in fact even
without this hypothesis, see 3.6, Step 3).}

In particular, if\/ $X$ satisfies $S_2$ and $\w_X$ is invertible then
 $$
\pi_*\w_Y=\Cee\cdot\w_X\qquad\text{and}\qquad\pi^*\w_X=\CHom(\Cee,\w_Y)=\w_Y(C),
 $$
where $\w_Y(C)$ is the $S_2$isation of\/ $\w_Y\otimes\Oh_Y(C)$, that is, the
divisorial sheaf on $Y$ of rational sections of $\w_Y$ with poles along $C$.
\endproclaim

\demo{{\rm2.4} Proof} The adjunction formula in the first sentence is standard
use of duality (see for example Proposition~2.11 of the Appendix for the case
of projective $k$@-schemes). A homomorphism $\al\:\pi_*\Oh_Y\to\w_X$ is of
course determined by $\al(1)$, so that
 $$
\CHom_{\Oh_X}(\pi_*\Oh_Y,\w_X)=
\bigl\{s\in\w_X\mid fs\in\w_X\text{\ for all\ } f\in
\pi_*\Oh_Y\bigr\},
 $$
which is equal to $\Cee\cdot\w_X$ on the locus where $\w_X$ is locally free.
The rest is easy: $\Cee$ satisfies $S_2$ by the above, and $\w_Y$ by
\cite{C3-f, App\. to \S1, Theorem~7}; and two coherent sheaves satisfying $S_2$
that coincide in codimension $1$ are equal. \QED\enddemo

\subheading{{\rm2.5}\enspace Etymology} The name {\it subadjunction} is
explained as follows: for an irreducible plane curve $X$ of degree $d$, the
canonical class of the resolution or normalisation $Y$ is the sheaf of
$\Oh_Y$@-modules generated by $\Cee\cdot\Oh(d-3)$, and the conductor ideal
$\Cee$ itself is determined in terms of {\it adjoint} forms, that is, forms
vanishing to order $m-1$ at every $m$@-fold point of $X$, including infinitely
near points (this goes back to Brill and Noether around 1870, and is also the
subject in the 1950s of Gorenstein's thesis and sections in Kodaira's papers
on surfaces). In higher dimension, the canonical class and plurigenera of the
resolution of an irreducible hypersurface with arbitrary singularities
correspond to {\it adjunction ideals} that can't be described in such simple
terms; however, by the proposition, the canonical class of the {\it
normalisation} is determined by the conductor ideal, given exactly as in the
curve case by conditions in codimension $1$, the {\it subadjunction}
conditions.

Subadjunction played a foundational role in the dark ages before K\"ahler
differentials and the Grothendieck dualising sheaf: the canonical class of a
nonsingular projective variety $V$ was often defined in terms of subadjunction
applied to a generic projection of $V$ as a hypersurface. It's clear from
Enri\-ques' discussion of subadjunction and adjunction for singular surfaces
in $3$@-space in \cite{Enri\-ques, Ch\.~III, \S\S6--7} that he understood
pretty well the case of ordinary multiple points, and the difficulties of
working with worse isolated singularities; but it's curious that he does not
seem to know Du~Val's work, the most substantial result known at the time.

 \myno{2.6} The next result solves the converse problem of determining
$\w_X$ in terms of $Y$ and $\fie\:C\to D$.

 \proclaim{Theorem} Assume that $X$ is $S_2$, so that in particular $C$ and $D$
have pure dimension $n-1$. Then applying the cohomological $\partial$@-functor
$\CExt_{\Oh_X}^*(\text{--},\w_X)$ to the exact diagram $(*)$ at the start of
\S2 gives the commutative diagram with exact rows
 $$
\matrix\format &\,\c\,\\
0&\to&\w_X&\to&\pi_*\w_Y(C)
&\to&\CExt^1_{\Oh_X}(\Oh_D,\w_X)&\to&0\\\vspace{4pt}
&&\bigcup&&\Vert&&\phantom{\scriptstyle{\Tr}\kern-2em}
\uparrow\scriptstyle{\Tr}\kern-2em\\\vspace{-.2pt}
0&\to&\pi_*\w_Y&\to&\pi_*\w_Y(C)&@>{\Res}>>
&\CExt^1_{\Oh_X}(\fie_*\Oh_C,\w_X)&\to&0
\endmatrix
\tag{$**$}
 $$
Moreover,
 $$
\align
\w_D&={}\text{$S_2$isation of\/ }\CExt^1_{\Oh_X}(\Oh_D,\w_X)\\
\fie_*\w_C&={}\text{$S_2$isation of\/ }\CExt^1_{\Oh_X}(\fie_*\Oh_C,\w_X);
\endalign
 $$
the $S_2$isation does nothing if $X$ is Cohen--Macaulay, so that
 $$
\w_D=\CExt^1_{\Oh_X}(\Oh_D,\w_X)
\quad\text{and}\quad
\fie_*\w_C=\CExt^1_{\Oh_X}(\fie_*\Oh_C,\w_X).
 $$

In more detail, this means the following:

(1) The composite of the surjection in the exact sequence
 $$
0\to\pi_*\w_Y\to\pi_*\w_Y(C)\to\CExt^1_{\Oh_X}(\fie_*\Oh_C,\w_X)\to0
 $$
and the $S_2$isation $\CExt^1_{\Oh_X}(\fie_*\Oh_C,\w_X)\subset\fie_*\w_C$ is a
canonically defined ``Poincar\'e residue'' map
$\Res=\Res_{Y,C}\:\w_Y(C)\to\w_C$.

(2) The trace map $\Tr=\Tr_{C/D}\:\fie_*\w_C\to\w_D$, which is
canonically defined and surjective, fits into a diagram
 $$
\spreadmatrixlines{4pt}
\matrix
\CExt^1_{\Oh_X}(\Oh_D,\w_X)&\subset&\w_D\\
\phantom{\scriptstyle{\Tr}\kern-2em}
\uparrow\scriptstyle{\Tr}\kern-2em
&&
\phantom{\scriptstyle{\Tr}\kern-2em}
\uparrow\scriptstyle{\Tr}\kern-2em\\
\CExt^1_{\Oh_X}(\fie_*\Oh_C,\w_X)&\subset&\fie_*\w_C.
\endmatrix
 $$

$(3)$ In these terms,
 $$
\w_X=\ker\left\{\pi_*\w_Y(C)@>{\Res}>>\fie_*\w_C@>{\Tr}>>\w_D
\right\}.
 $$\endproclaim

\demo{{\rm2.7} Proof} This all follows formally from the adjunction
properties of predualising sheaves (discussed in the appendix). First of all,
clearly $\CHom(\Oh_D,\w_X)=\CHom(\Oh_C,\w_X)=0$ since $\w_X$ is torsion-free;
$\CHom_{\Oh_X}(\Oh_X,\w_X)=\w_X$ is obvious, and subadjunction gives
$\CHom_{\Oh_X}(\pi_*\Oh_Y,\w_X)=\pi_*\w_Y$. Next,
 $$
\CHom_{\Oh_X}(\Cee,\w_X)=\CHom_{\Oh_X}(\pi_*\Cee_Y,\w_X)=\pi_*\w_Y(C).
 $$
This holds because $\Cee=\pi_*\Cee_Y$ is a $\pi_*\Oh_Y$@-module, so any
$\Oh_X$@-linear homo\-morphism from it can only map into a
$\pi_*\Oh_Y$@-submodule
of $\w_X$, that is, into $\pi_*\w_Y$; thus the middle term is
 $$
\CHom_{\pi_*\Oh_Y}(\pi_*\Cee_Y,\pi_*\w_Y)=\pi_*\w_Y(C).
 $$

Next, $\CExt^1_{\Oh_X}(\Oh_X,\w_X)=0$ follows from basic properties
of $\Ext$s, since the stalks of $\Oh_X$ are projective. As just explained, the
two functors $\CHom_{\Oh_X}(\text{--},\w_X)$ and
$\CHom_{\pi_*\Oh_Y}(\text{--},\pi_*\w_Y)$ coincide on $\pi_*\Oh_Y$@-modules,
so that the fact that the stalks of $\pi_*\Oh_Y$ are projective over
$\pi_*\Oh_Y$ implies that also $\CExt^1_{\Oh_X}(\pi_*\Oh_Y,\w_X)=0$.

Finally, the adjunction formula Theorem~2.12 of the Appendix gives that $\w_D$
is the $S_2$isation of $\CExt_{\Oh_X}^1(\Oh_D,\w_X)$ and $\fie_*\w_C$ that of
$\CExt_{\Oh_X}^1(\fie_*\Oh_C,\w_X)$.

Therefore, applying $\CExt_{\Oh_X}^*(\text{--},\w_X)$ to $(*)$ gives the
diagram $(**)$ of the theorem. \QED\enddemo

\proclaim{{\rm2.8}\enspace Corollary} Equivalent conditions:

(i) The predualising sheaf $\w_X$ is an invertible $\Oh_X$@-module.

(ii) $\w_X$ is invertible in codimension $1$, and for all $P\in Y$ there
exists an element $s\in\ker(\Tr_{C/D}\circ\Res_{Y,C})\:\pi_*\w_Y(C)\to\w_D$
such that $\w_Y(C)=\Oh_Y\cdot s$ near $P$.

(iii) $\w_X$ is invertible in codimension $1$, $\w_Y(C)$ is invertible, and
for all $P\in C$ there exists a basis element $s\in\ker\Tr_{C/D}\:\allowbreak
\CExt^1_{\Oh_X}(\fie_*\Oh_C,\w_X)\to\CExt^1_{\Oh_X}(\Oh_D,\w_X)$ such that
$\CExt^1_{\Oh_X}(\fie_*\Oh_C,\w_X)=\Oh_C\cdot s$ near $P$.

(iv) $\w_Y(C)$ is an invertible $\Oh_Y$@-module and
 $$
\ker\{\Tr_{C/D}\:\allowbreak
\CExt^1_{\Oh_X}(\fie_*\Oh_C,\w_X)\to\CExt^1_{\Oh_X}(\Oh_D,\w_X)\}
$$
is an invertible $\Oh_D$@-module.

\rm Note that when $X$ is Cohen--Macaulay, conditions (iii) and (iv) refer
simply to the kernel of $\Tr_{C/D}\:\fie_*\w_C\to\w_D$
\endproclaim

\demo{Proof} (i)$\implies$(ii) or (iv) is clear, since $\w_X=\Oh_X\cdot s$
implies that $\w_Y(C)=\pi^*\w_X=\Oh_Y\cdot s$ by Proposition~2.3, and
$s\in\ker(\Tr\circ\Res)$ by Theorem~2.6; and $s$ maps to an $\Oh_D$@-basis of
$\ker\Tr_{C/D}$. (ii)$\iff$(iii) is clear by the surjectivity of $\Res$ in
$(**)$, since if $s\in\w_Y(C)$ maps to $\overline s=\Res
s\in\CExt^1_{\Oh_X}(\fie_*\Oh_C,\w_X)$ then $s$ is an $\Oh_Y$@-basis of
$\w_Y(C)$ if and only if $\overline s$ is an $\Oh_C$@-basis of
$\CExt^1_{\Oh_X}(\fie_*\Oh_C,\w_X)$.

(ii)$\implies$(i) Let $Q\in X$. First suppose that
$s\in\ker(\Tr_{C/D}\circ\Res_{Y,C})$ is an $\Oh_Y$@-basis of $\w_Y(C)$ near $P$
for every $P$ lying over $Q$; required to prove that $\Oh_X\cdot s=\w_X$ near
$Q$. This is true after localising at the generic point of a codimension $1$
subvariety $\Ga\subset X$: for by assumption, $\w_{X,\Ga}$ is locally free,
and if $s\in m_\Ga\cdot\w_{X,\Ga}$ then $s$ could not be a basis of $\w_Y(C)$
above $\Ga$. The result then follows from the fact that $\w_X$ is $S_2$.

By (ii), for each $P\mapsto Q$ there is an element
 $$
s_P\in\ker(\Tr_{C/D}\circ\Res_{Y,C})\:\pi_*\w_Y(C)\to\w_X
$$
that is an $\Oh_Y$@-basis of $\w_Y(C)$ near $P$. If the residue field of
$\Oh_{X,Q}$ is infinite, or big enough compared with the number of $P$, then
a suitable linear combination of the $s_P$ with coefficients in $\Oh_{X,Q}$
will be a basis at every $P$.

It's not hard to deal with the case that $\Oh_{X,Q}$ has a finite residue
field by making a finite faithfully flat extension
$\Oh_{X,Q}\subset\Oh_{X',Q}$ to increase the residue field, concluding that
$\w_{X'}$ is locally free by the above argument, and then using the fact that
a finite module over a local ring which becomes free after a faithfully flat
extension was already free (because free = flat, e.g., \cite{Bourbaki, Cor\.
II.3.2.5.2}).

(iv)$\implies$(i) Outside $D$ and $C$, obviously $\w_X=\w_Y=\w_Y(C)$; near
$D$, the exact diagram $(**)$ can be rewritten
 $$
\spreadmatrixlines{3pt}
\matrix\format &\,\c\,\\
0&\to&\pi_*\w_Y&\to&\w_X&\to&\ker\Tr_{C/D}&\to&0\\
&&\Vert&&\bigcap&&\bigcap\\
0&\to&\pi_*\w_Y&\to&\pi_*\w_Y(C)&\to
&\CExt^1_{\Oh_X}(\fie_*\Oh_C,\w_X)&\to&0.\\
\endmatrix
 $$
The kernel in the top left is $\w_X\cap\pi_*\w_Y=\pi_*\w_Y$.

Pick an $\Oh_D$@-basis $\overline s\in\ker\Tr_{C/D}$, and $s\in\w_X$ mapping to
it. I claim that $\Oh_X\cdot s=\w_X$. Since both $\Oh_X\cdot s$ and $\w_X$ are
$S_2$, it's enough to prove this in codimension $1$. Now $\Cee\cdot\w_X\subset
\pi_*\w_Y$ by Proposition~2.3, and I will prove in 3.6, Step~3 that equality
holds in codimension $1$. Now on the locus where $\pi_*\w_Y=\Cee\cdot\w_X$, I
have
$\ker\Tr_{C/D}=\w_X/(\Cee\cdot\w_X)=\w_X\otimes\Oh_D$, so that
 $$
\overline s\;\text{bases }\ker\Tr_{C/D}\implies s\;\text{bases}\;\w_X
 $$
follows from Nakayama's lemma. \QED\enddemo

\subheading{{\rm2.9}\enspace Remark} Writing down the trace map $\Tr_{C/D}$ is
an
activity that takes place at the generic point of each component of $D$, and
reduces there to duality for Artinian local rings. Each of $\fie_*\w_C$ and
$\w_D$ is torsion-free, so a subsheaf of its rational sections (the sheaf made
up of direct sums of generic stalks). If $\Ga$ is a component of $D$, then
the generic stalk $\Oh_{D,\Ga}$ is an Artinian local ring, with dualising
module $\w_{D,\Ga}$, and $\Tr\:\fie_*\w_{C,\Ga}\to\w_{D,\Ga}$ is just the
$\Oh_{D,\Ga}$@-dual of the inclusion $\Oh_{D,\Ga}\subset\fie_*\Oh_{C,\Ga}$.
Thus $\ker\Tr_{C/D}$ can be described as the subsheaf
 $$
\ker\Tr_{C/D}=\{s\in\fie_*\w_C\mid\Tr(s_\Ga)=0\in\w_{D,\Ga}\text{ for all
}\Ga\}\subset\fie_*\w_C.
 $$

Although $\Tr_{C/D}$ has this ``birational'' description, the condition for
$\ker\Tr_{C/D}$ to be invertible as an $\Oh_D$@-module is nevertheless a
delicate biregular question, and understanding this in a special case is the
main point of \S4. A subtle point that causes a lot of confusion is that
$\Tr_{C/D}$ is only linear over $\Oh_D$; however, the subring $\Oh_D$ is the
{\it unknown} in my calculations in \S4.

\heading Appendix to \S2. Adjunction for a finite morphism\endheading

\myno{2.10} This section is a technical digression. $\pi\:Y\to X$ is a finite
morphism of schemes, not necessarily surjective. By \cite{Hartshorne, Ch\.~II,
Ex\. 5.17, (e)}, $\pi_*$ identifies the category of $\Oh_Y$@-modules on
$Y$ with that of $\Oh_X$@-modules on $X$ together with an action of
$\pi_*\Oh_Y$; the map back is the module-to-sheaf construction
$M\mapsto\widetilde M$ generalised from affine schemes to affine morphisms. I
just say $\Oh_Y$@-module from now on.

The right adjoint of $\pi_*$ takes an $\Oh_X$@-module $\G$ into
$\CHom_{\Oh_X}(\pi_*\Oh_Y,\G)$, made into an $\Oh_Y$@-module by multiplication
in the source by elements of $\pi_*\Oh_Y$: as written in Holy Scriptures
\cite{Grothendieck--Hakim, Grothendieck--Hartshorne, \S4}, if $\Cee$ is the
category of $A$@-modules, and $T\:\Cee\to\roman{\underline{Ab}}$ a
contravariant functor, then $T(M)$ has a canonical $A$@-module structure with
$a\in A$ acting by $T(\mu_a)$, where $\mu_a$ is the homo\-thety of
multiplication by $a$; this must be borne in mind throughout. For an
$\Oh_Y$@-module $\F$ and an $\Oh_X$@-module $\G$, there is a canonical
bifunctorial isomorphism
 $$
\Hom_{\Oh_X}(\pi_*\F,\G)=\Hom_{\Oh_Y}(\F,\CHom_{\Oh_X}(\pi_*\Oh_Y,\G)).
 $$

To use easy characterisations of the predualising sheaf $\w_X$ I work here
only with quasiprojective $k$@-schemes. A more sophisticated definition of the
dualising complex makes everything work more generally, somewhat
tautologically; but the existence of the dualising complex is hard in general.

Recall from \cite{Hartshorne, Ch\.~III, \S7} that the predualising sheaf
$\w_X$ on an $n$@-dimensional projective scheme $X$ is determined by the
following universal mapping property: there is a $k$@-linear map ``trace''
$t\:H^n(X,\w_X)\to k$, and for every coherent sheaf $\F$ on $X$, any
$k$@-linear map $H^n(X,\F)\to k$ is induced by a morphism $\F\to\w_X$, so that
$H^n(X,\F)\dual\Hom(\F,\w_X)$.

\proclaim\nofrills{{\rm2.11}\enspace Proposition}{ \rm (\cite{Hartshorne,
Ch\.~III, Ex.~7.2}).} Let $\pi\:Y\to X$ be a finite morphism of projective
$k$@-schemes with $\dim X=\dim Y$. Then
 $$
\w_Y=\CHom_{\Oh_X}(\pi_*\Oh_Y,\w_X)
$$
and the map $\w_Y\to\w_X$ defined by $\alpha\mapsto\alpha(1)$ is the trace map
for $\pi$.
\endproclaim

\demo{Proof} Let $\F$ be any coherent sheaf on $Y$; then (omitting the $\pi_*$
for clarity), I get
 $$
\align
H^n(Y,\F)=H^n(X,\F)\dual{}&\Hom_{\Oh_X}(\F,\w_X)\\
={}&\Hom_{\Oh_Y}(\F,\CHom_{\Oh_X}(\Oh_Y,\w_X).
\endalign
 $$
Therefore $\w_Y=\CHom_{\Oh_X}(\Oh_Y,\w_X)$ satisfies the universal
property of a dualising sheaf on $Y$.\QED\enddemo

\proclaim\nofrills{{\rm2.12}\enspace Theorem}{\rm \; (Adjunction formula).\ }\
Let $X$ be a purely $n$@-dimensional quasiprojective scheme, $\w_X$ a
predualising sheaf for $X$ \cite{Hartshorne, Ch\.~III, \S7} and $\pi\:Y\to X$
a finite morphism, with $Y$ purely of dimension $n-r$.

(1) Suppose that $X$ is Cohen--Macaulay. Then the predualising sheaf of $Y$
is given by
 $$
\w_Y=\CExt_{\Oh_X}^r(\pi_*\Oh_Y,\w_X).
 $$

(2) Suppose that $X$ is Cohen--Macaulay at every codimension $1$ point of $Y$
(this holds in particular if $X$ satisfies Serre's condition $S_{r+1}$). Then
$\w_Y$ is the $S_2$isation (see 2.2) of $\CExt_{\Oh_X}^r(\pi_*\Oh_Y,\w_X)$.

\endproclaim

\demo{Proof} (1) Suppose $X\subset\proj^N=\proj$ has codimension $s$. Then
 $$
\w_Y=\CExt_{\Oh_\proj}^{r+s}(\Oh_Y,\w_\proj).
 $$
This $\CExt$ is the $(r+s)$th homology sheaf of the complex
$\CHom_{\Oh_\proj}(\Oh_Y,\I\hidot)$, where $\w_\proj\to\I\hidot$ is an
injective resolution. But $\Oh_Y$ is a sheaf of $\Oh_X$@-modules, so that by
the above, I can $\CHom$ via $\Oh_X$:
 $$
\CHom_{\Oh_\proj}(\Oh_Y,\I\hidot)=
\CHom_{\Oh_X}\bigl(\Oh_Y,\CHom_{\Oh_\proj}(\Oh_X,\I\hidot)\bigr).
 $$
I claim that the inner complex $\J\hidot=\CHom_{\Oh_\proj}(\Oh_X,\I\hidot)$ is
essentially an injective resolution of $\w_X$ shifted by $s$. Indeed, it is
easy to see that $\I^i$ injective over $\Oh_\proj$ implies that
$\J^i=\CHom_{\Oh_\proj}(\Oh_X,\I^i)$ is injective as $\Oh_X$ module. Also,
because $X$ is Cohen--Macaulay, the complex $\J\hidot$ has cohomology
 $$
\CExt_{\Oh_\proj}^i(\Oh_X,\w_\proj)=\cases
\w_X&\text{if }i=s,\\
0& \text{if }i\ne s.\\
\endcases
 $$
Thus $\J\hidot$ can be written as a direct sum of two complexes
$\J\hidot_1\oplus\J\hidot_2$, where $\J\hidot_1$ is an exact complex of
injective $\Oh_X$@-modules in degrees $[0,s]$ and $\J\hidot_2$ an injective
resolution of $\w_X$ starting in degree $s$. Thus
 $$
\CExt^{r+s}_{\Oh_\proj}(\Oh_Y,\w_\proj)=\CExt^r_{\Oh_X}(\Oh_Y,\w_X).
$$

(2) follows at once from the fact that $\w_Y$ is $S_2$, and by what I've
just said, coincides with $\CExt^r_{\Oh_X}(\Oh_Y,\w_X)$ in codimension 1.
\QED\enddemo

\subheading{{\rm2.13}\enspace Examples} Here is the well-known case of the
theorem with $r=1$, and a counterexample pointed out by Koll\'ar to show that
the $S_2$isation is really needed in Theorem~2.12.

\smallskip
(a) Suppose $Y\subset X$ is a Cartier divisor; then applying $\CExt$s to the
exact sequence $0\to\Oh_X(-Y)\to\Oh_X\to\Oh_Y\to0$, gives
 $$
0\to\CHom(\Oh_X,\w_X)\to\CHom(\Oh_X(-Y),\w_X)\to
\CExt_{\Oh_X}^1(\Oh_Y,\w_X)\to0,
 $$
that is, $\CExt_{\Oh_X}^1(\Oh_Y,\w_X)=\w_X(Y)/\w_X=\w_X(Y)\otimes\Oh_Y$. Hence
if $X$ is $S_2$ then
 $$
\w_Y=\text{$S_2$isation of }\CExt_{\Oh_X}^1(\Oh_Y,\w_X)
=\text{$S_2$isation of }\w_X(Y)\otimes\Oh_Y.
 $$

\smallskip
(b) Let $X$ be the projective cone over a normally embedded Abelian surface
$A$ and $Y$ a hypersection through the vertex point $O$. Then
$\w_X\cong\Oh_X(-Y)$ and $\CExt_{\Oh_X}^1(\Oh_Y,\w_X)\cong\Oh_Y$ by the
argument of (a) above. But $Y$ is not normal, so that $\Oh_Y$ is not $S_2$;
here $\w_Y=\Oh_{\wave Y}$ where $\pi\:\wave Y\to Y$ is the normalisation.
\smallskip

(c) Let $Y=Y_1+Y_2$ be Cartier divisors on $X$ as in (a). Then
the theorem applies to the inclusion morphism $Y_1\subset Y$, giving
 $$
\w_{Y_1}=\CHom_{\Oh_Y}(\Oh_{Y_1},\w_Y).
 $$
It's not hard to evaluate this to be $=\w_Y(-Y_2)\otimes\Oh_{Y_1}$.

\heading\S3. What happens in codimension $1$\endheading

Duality for finite modules over a local Artinian (0@-dimensional) ring is
almost
as easy and explicit as for vector spaces over a field. The aim of this
section is to apply this duality to complete the proof of Corollary~2.8, (ii)
or (iv)$\implies$(i), to translate the condition for $X$ to be Gorenstein in
codimension $1$ into a more explicit form, and to solve this condition in the
simplest cases.

Most of the material is copied more or less directly from \cite{Serre, Ch\.~IV,
\S11}, although my category is closer to the natural level of generality sought
by the true Bourbakist.

\subheading{{\rm3.1}\enspace Notation and assumptions} Except where stated
otherwise,
I localise everything throughout \S3 at a codimension $1$ point of $X$, that
is, at the generic point of an irreducible codimension $1$ locus
$\Ga\subset X$. Suppose that $C\subset Y$ localised above $\Ga\subset X$
has the decomposition into prime divisors $C=\sum n_E E$. Thus $\Oh_X$ is a
reduced $1$@-dimensional local ring, and is a subring of $\Oh_Y$, which is a
product of DVRs:
 $$
\Oh_Y=\prod R_E,\quad\text{where $R_E=\Oh_{Y,E}$ is a DVR with local parameter
$t_E$}.
 $$
I suppress $\pi_*$ throughout \S3. Write $m=m_\Ga\subset\Oh_X$ for the
maximal ideal and $K=k(\Ga)=\Oh_X/m$ for its residue field, the function
field of the subscheme $\Ga$; note that the residue field $k(E)=R_E/(t_E)$
of each localisation $R_E$ of $\Oh_Y$ is a finite extension of $K$. The
geometric
picture of the normalisation is ``unzipping'' $X$ along its codimension $1$
singular locus, and in general the set-theoretic multiple locus of $Y\to X$ can
be a ramified (or even inseparable) cover $\coprod E\to\Ga$.

The conductor $\Cee$ has finite colength in $\Oh_X$ and $\Oh_Y$, and the
quotients $\Oh_D=\Oh_X/\Cee$ and $\Oh_C=\Oh_Y/\Cee$ are Artinian rings, with
$\Oh_D$ local and $\Oh_C=\prod A_E$ a product of rings of the form
 $$
A_E=\Oh_{nE}=R_E/(t_E{}^{n_E})\qquad\text{with $n_E\ge1$}.
 $$
Since $\Oh_D$ is Artinian and local, the stalk $\w_D$ of the dualising sheaf at
the generic point of $\Ga$ is the dualising module of $\Oh_D$, that is,
$\w_D=\Hom_{k_0}(\Oh_D,k_0)$ if $\Oh_D$ contains a field $k_0$ such that
$k_0\subset K$ if a finite extension (hands up those who don't remember how the
$\Oh_D$@-module structure in the $\Hom$ is defined!); and, quite generally,
$\w_D$ is the injective hull of the residue field $K$ as $\Oh_D$@-module. Thus
$\Hom_{k_0}(\text{--},k_0)$ in the $k_0$@-algebra case, or
$\Hom_{\Oh_D}(\text{--},\w_D)$ quite generally, is the dualising functor for
$\Oh_D$@-modules, and I write $M\dual N$ to mean $M=\Hom_{\Oh_D}(N,\w_D)$.

In this section the {\it length} $\ell(M)$ of a module $M$ {\it always means
its length as an $\Oh_D$@-module}, so that, for example, the residue field
$k(R_E)$ has length $\ell(k(R_E))=[k(R_E):K]$, and $A_E$ in the preceding
display has length $\ell(A_E)=n_E\cdot[k(R_E):K]$; as an $A_E$@-module, this
has of course a Jordan--H\"older composition series of length only $n_E$
defined by the powers of $t_E$, with successive quotients $k(R_E)$.

\proclaim{{\rm3.2}\enspace Main theorem} In the notation and assumptions of
\S2, consider the two diagrams
 $$
\spreadmatrixlines{.3cm}
\matrix\format\r\,&\,\c\,&\,\c\,&\,\c\,&\,\l\\
\Oh_Y(-C)=\Cee&\subset&\Oh_X&\subset&\Oh_Y\\
\w_Y(C)&\supset&\w_X&\supset&\w_Y\\
\endmatrix
\qquad\text{and}\qquad
\matrix\format\r\,&\,\c\,&\,\c\,&\,\c\,&\,\l\\
0&\subset&\Oh_D&\subset&\Oh_C\\
\w_C&\supset&\ker\Tr_{C/D}&\supset&0,\\
\endmatrix
 $$
the first of which consists of sheaves of $\Oh_X$@-modules, and the second of
$\Oh_D$@-modules, the corresponding quotients by $\Cee$ or $\w_Y$.

(I) Localised in codimension $1$, the aligned
$\subset$s and $\supset$s correspond to $\Oh_D$@-dual modules
 $$
\align
\Oh_X/\Cee=\Oh_D&\dual\w_Y(C)/\w_X=\w_D,\\ \vspace{.2cm}
\Oh_Y/\Oh_X=\Oh_C/\Oh_D&\dual\w_X/\w_Y=\ker\Tr_{C/D}.
\endalign
 $$

(II) The length of the $\Oh_D$@-modules in the top line is $\le$ that in
the bottom line, $\ell(\Oh_D)\le\ell(\Oh_C/\Oh_D)=\ell(\w_X/\w_Y)$, and
 $$
\text{equality}\;\iff\;\Oh_D\cong\ker\Tr_{C/D}=\w_X/\w_Y\;\iff\;\Oh_X\cong\w_X,
 $$
that is, if and only if\/ $\w_X$ is Gorenstein in codimension $1$.

{\rm In general (not localising at a codimension $1$ point), the first diagram
can be taken to mean duality of reflexive sheaves.}\endproclaim

\subheading{{\rm3.3}\enspace Remark} In (I), all the equalities are already
known
from \S2, with those on the left by definition, and those on the right coming
from Theorem~2.6, so that the point is to prove the duality. (II) is the famous
inequality $n\le2\de$ of \cite{Serre, Ch\.~IV, \S11}, where
 $$
n=\ell(\Oh_Y(C)/\Oh_Y)=\ell(\Oh_C/\Oh_D)+\ell(\Oh_D)=\sum n_E\cdot [k(E):K]
 $$
(here $C=\sum n_E E$, with the $E$ prime divisors of $Y$), and
$\de=\ell(\w_X/\w_Y)$; if $\pi\:Y\to X$ is the normalisation (resolution) of
a singularity $P\in X$ of an irreducible projective curve over an algebraically
closed field, then $\de$ is the {\it genus} of the singularity $P\in X$, that
is,
 $$
\de=\de p_a=h^0(\w_X)-h^0(\w_Y).
 $$

\myno{3.4} The proof of the theorem follows closely the arguments of
\cite{Serre, Ch\.~IV, \S11}, which are mainly in terms of $\pi\:Y\to X$, but
it's convenient to waste a page of journal space discussing the following
definition, which abstracts away from the $1$@-dimensional schemes $X$ and $Y$,
and focuses on the Artinian conductor subschemes $D$ and $C$ where the magic
duality works.

\subheading{Definition \rm (half-filling of $\Oh_C$)} Suppose given a field
$K$, and a number of DVRs $R_E$, with local parameters $t_E$, whose residue
fields $k(E)=R_E/(t_E)$ are all finite extensions of $K$; let $\Oh_C=\prod A_E$
be a product of rings of the form
 $$
A_E=R_E/(t_E{}^{n_E})\qquad\text{with $n_E\ge1$}.
 $$

A {\it part-filling} of $\Oh_C$ is a subring $\Oh_D\subset\Oh_C$ satisfying

(i) $\Oh_D$ is a local ring with residue field $\Oh_D/m_D=K$;

(ii) $\Oh_C/\Oh_D$ is a faithful $\Oh_D$@-module, that is,
$\Ann_{\Oh_D}(\Oh_C/\Oh_D)=0$.

\noindent If in addition

(iii) $\Oh_C/\Oh_D\cong\w_D$ as $\Oh_D$@-module,

\noindent then $\Oh_D$ is a {\it half-filling of} $\Oh_C$.

\subheading{{\rm3.5}\enspace Remarks} (a) In practical calculations, the
trickiest thing to work with is the condition that $\Oh_D$ is a subring, since
this is nonlinear; see Proposition~3.9 and 4.4 for an example.

(b) The local condition (i) means that the image of $\Oh_D$ under the
composite map $\Oh_C\to\prod R_E\to\prod K_E$ is the diagonal $K\subset\prod
K_E$. This is equivalent to saying that $\dim\Hom_{\Oh_D}(\Oh_D,K)=1$. If $K$
is algebraically closed, it's equivalent to any of the following: the identity
element $(1,\dots,1)\in\Oh_C$ is in $\Oh_D$, but no other idempotent; every
element $(f_1,\dots,f_r)\in m_D$ has each $f_E$ in the maximal ideal $t_EA_E$
of $A_E$; every nonunit of $\Oh_D$ is nilpotent.

(c) If $Y\to X$ is a normalisation localised in codimension $1$ as at the start
of \S3, the conductor subschemes $C\subset Y$ and $D\subset X$ define a
part-filling $\Oh_D\subset\Oh_C$, with (ii) the definition of the conductor;
(iii) will be equivalent to $X$ Gorenstein. (iii) obviously implies (ii) and
also $\ell(\Oh_D)=\ell(\Oh_C/\Oh_D)$. The main point to prove in Theorem~3.2 is
that
 $$
\text{(ii)}\implies\ell(\Oh_D)\le\ell(\Oh_C/\Oh_D),
\qquad\text{and}\qquad\text{equality}\implies\text{(iii)}.
 $$

(d) Every part-filling arises as a conductor subscheme for a normalisation
$\pi\:Y\to X$ of a reduced Noetherian local scheme $X$: just take $X$ and $Y$
to be $\Spec$ of
 $$
\Oh_Y=\prod R_E\quad\text{and}\quad\Oh_X=\{f\in\Oh_Y\mid f\mapsto\bar
f\in\Oh_D\subset\Oh_C\}.
 $$
In the purely geometric case, when $\Oh_D$ is a $K$@-algebra and $K=K_E$, you
can
even assume that each component of $Y$ is $\cong\aff^1_K$, but the more general
category also allows amusing things like reduced divisors on an arithmetic
surface with some components in fibres (characteristic $p$) and some horizontal
(mixed characteristic). So scheme theory suggests, for example, glueing
together
a ring of integers in a number field and a projective curve over a finite
field,
constructions which a number theorist may not immediately think of as natural.

(e) I've forgotten what I wanted to say here. Oh yes, the definition could in
principle be generalised by allowing $K$ and the $K_E$ to be finite extensions
of
a common subfield $k_0$, but the extra generality is illusory, e.g., by
\cite{Matsumura, Theorem~28.3, (ii)} in the equal characteristic case.

(f) At cherry blossom time in Kyoto, in connection with the inequality in
(II), I conjectured foolishly that an Artinian ring $A$ having a faithful
module $M$ of finite length should satisfy $\ell(A)\le\ell(M)$. A
counterexample: the ring of $2n\times2n$ matrixes having zero entries except
in the top right $n\times n$ block, and equal diagonal entries, is clearly a
commutative local $k$@-algebra of length $n^2+1$ having $k^{2n}$ as a faithful
module.

\myno{3.6} I now proceed to prove Theorem~3.2 and other matters of interest in
5 easy steps. Each step consists of a statement for an abstract part-filling
$\Oh_D\subset\Oh_C$, that makes sense without mention of $X$ and $Y$, then
more-or-less trivial consequences for $Y\to X$.

\proclaim{Step~1} $\ker\Tr_{C/D}$ is the $\Oh_D$@-dual of $\Oh_C/\Oh_D$,
that is,
 $$
\Oh_C/\Oh_D\dual\ker\{\Tr_{C/D}\:\w_C\to\w_D\}=\w_X/\w_Y.
 $$\endproclaim

\demo{Proof} In the case when $\Oh_C$ and $\Oh_D$ are finite $k_0$@-algebras,
their dualising modules are $\w_C=\Hom_{k_0}(\Oh_C,k_0)$ and
$\w_D=\Hom_{k_0}(\Oh_D,k_0)$, and the trace map
 $$
\Tr\:\Hom_{k_0}(\Oh_C,k_0)\to\Hom_{k_0}(\Oh_D,k_0)
 $$
is just the restriction from $\Oh_C$ to $\Oh_D$, that is, the dual of the
inclusion map $\Oh_D\hookrightarrow\Oh_C$, so clearly
 $$
\ker\Tr=\Hom_{k_0}(\Oh_C/\Oh_D,k_0).
 $$
The general case is exactly the same on replacing $\Hom_{k_0}(\text{--},k_0)$
by $\Hom_{\Oh_D}(\text{--},\w_D)$ throughout. \QED\enddemo

\proclaim{Step~2} Therefore
$\ell(\Oh_C/\Oh_D)=\ell(\ker\Tr_{C/D})=\ell(\w_X/\w_Y)$ and the $\Oh_D$@-module
$\ker\Tr_{C/D}$ is faithful, that is, $\Ann_{\Oh_D}(\ker\Tr_{C/D})=0$;
hence $$
\Cee=\Ann_{\Oh_X}(\Oh_Y/\Oh_X)=\Ann_{\Oh_X}(\w_X/\w_Y).
 $$\endproclaim

\demo{Proof} This follows from Step~1, since the dual $\Oh_D$@-modules
$\Oh_C/\Oh_D$ and $\ker\Tr_{C/D}$ have the same length and the same
annihilator;
and the $\Oh_X$@-module structures come from $\Oh_X\twoheadrightarrow\Oh_D$.
\QED\enddemo

\proclaim{Step~3} For each $E$, write
$t'_E=(1,\dots,t_E,\dots,1)\in\Oh_C$; this corresponds to a function on $Y$
with zero of order $1$ along $E$, and no zeros along other components. Then for
each $E$, I claim that
 $$
\ker\Tr_{C/D}\not\subset t'_E\cdot\w_C,
 $$
that is, $\ker\Tr_{C/D}$ contains a basis of the localisation
$\w_{C,E}\cong\Oh_{C,E}=A_E$.

Thus $\w_X\subset\w_Y(C)$ contains a local basis of $\w_Y(C)$ at the generic
point of each $E$, that is, an element of $\w_Y(C)$ with pole of order exactly
$n_E$. Therefore $\w_Y(C)$ is the $S_2$isation of $\Oh_Y\cdot\w_X$; and
$\w_Y$ the $S_2$isation of $\Cee\cdot\w_X$, so the hypothesis on $X$ can be
omitted in the conclusion of Proposition~2.3. {\rm(Note that the last two
statements
are for a general variety before localising, but only the local statement in
codimension $1$ requires proof.)}\endproclaim

\demo{Proof} Consider the element $s'_E=(0,\dots,t_E^{n_E-1},\dots,0)\in\Oh_C$.
Assuming by contradiction that $\ker\Tr_{C/D}\subset t'_E\cdot\w_C$, I prove
that
 $$
0\ne s'_E\in\Ann_{\Oh_D}(\ker\Tr_{C/D}),
 $$
which contradicts Step~2. To do this, I must prove both that multiplication by
$s'_E$ is zero, and that $s'_E\in\Oh_D$. Obviously $0\ne s'_E\in\Oh_C$, and
the multiplication map $\mu_{s'_E}\:\Oh_C\to\Oh_C$ is $\Oh_C$@-linear,
therefore $\Oh_D$@-linear. The dual $\Oh_D$@-linear map $\w_C\to\w_C$ induced
by $\mu_{s'_E}$ kills $t'_E\cdot\w_C$, and therefore by the assumption also
$\ker\Tr_{C/D}=\Hom_{\Oh_D}(\Oh_C/\Oh_D,\w_D)$. Finally, by duality it follows
that $\mu_{s'_E}(\Oh_C)=\Oh_C\cdot s'_E\subset\Oh_D$, so that in particular
$s'_E\in\Oh_D$. This is the required contradiction.

The remaining assertions follow easily. \QED\enddemo

\proclaim{Step~4} $\ell(\Oh_D)\le\ell(\ker\Tr_{C/D})=\ell(\Oh_C/\Oh_D)=
\ell(\w_X/\w_Y)$; and if the residue field has cardinality
 $$
\card K>\text{\rm{number of components $E$ of $C$,}}
 $$
in particular if $K$ is infinite, then there exists an inclusion
$\Oh_D\hookrightarrow\ker\Tr_{C/D}=\w_X/\w_Y$.\endproclaim

\demo{Proof} For each $E$, I can pick an element $s_E\in\ker\Tr_{C/D}$ which is
a local basis of $\w_C$ at the generic point of $E$ by Step~3. Under the given
condition on the cardinality of $K$, it's easy to see that a suitable linear
combination $s=\sum f_Es_E$ with $f_E\in\Oh_D$ is a local basis of $\w_C$ at
the
generic point of every $E$; therefore $f\mapsto f\cdot s\in\w_C$ is an
injective
map $\Oh_C\to\w_C$ (in fact an isomorphism), so its restriction to $\Oh_D$ is
also injective.

This proves the inequality $\ell(\Oh_D)\le\ell(\Oh_C/\Oh_D)=\ell(\w_X/\w_Y)$ in
the case that the residue field $K$ is infinite. For small $K$, the inequality
(but not the inclusion) reduces without difficulty to the case of infinite $K$
by the standard trick \cite{Matsumura, p\. 114} of passing to the flat overring
$\Oh_D\mapsto\Oh_D(x)$, which makes the purely transcendental extension
$K\mapsto K(x)$ of the residue field. \QED\enddemo

\proclaim{Step~5} The equality
$\ell(\Oh_D)=\ell(\ker\Tr_{C/D})=\ell(\w_X/\w_Y)$ (or $n=2\de$ in
\cite{Serre, Ch\.~IV, \S11}) is a necessary and sufficient condition for
$\Oh_D$ to be a half-filling of $\Oh_C$, or for $X$ to be Gorenstein in
codimension $1$; dedicated algebraists can find a further halfdozen equivalent
conditions in terms of primary decomposition, socles and homological
algebra.\endproclaim

\demo{Proof} $\ker\Tr_{C/D}\dual\Oh_C/\Oh_D$ by Step~1, so (iii) in
Definition~3.4 implies that $\ell(\Oh_D)=\ell(\ker\Tr_{C/D})$. Conversely, if
the lengths are equal, the inclusion $\Oh_D\hookrightarrow\ker\Tr_{C/D}$ of
Step~4 must be an isomorphism, which proves the result if $K$ is big enough. As
before (Corollary~2.8, Proof of (ii)$\implies$(i)), to deal with the case of
finite residue field $K$ one has to make a finite faithfully flat extension
$\Oh_X$ inducing a sufficiently big residue field extension, and argue on a
finite module over a local ring which becomes free after a faithfully flat
extension. The last sentence is obvious.\QED\enddemo

\myno{3.7} The simplest application of these ideas is the description of
the codimension $1$ behaviour of nonnormal del Pezzo surfaces stated in the
introduction. Everything is still localised at a codimension $1$ point of $X$.

\proclaim{Theorem} (I) Let $\Oh_C=\prod R_E/(t_E{}^{n_E})$, and
assume that $n_E=1$ for some $E$. If $\Oh_C$ has a half-filling $\Oh_D$, then
$\Oh_D=k(\Ga)=K$ is a field (where $\Ga=D\red$), and\/ $\Oh_C$ is
either $\cong K\times K$, or a quadratic extension field of\/ $K$.

That is, if\/ $X$ is nonnormal and Gorenstein, and $C\subset Y$ has a reduced
component $E$, then $X$ has ordinary double points in codimension $1$ along
$\Ga$. {\rm In characteristic $2$, this includes {\it inseparable ordinary
double points}, that is, the quadratic extension may be an inseparable cover
$E\to\Ga$, in which case every geometric transverse section of the
singularity of $X$ along $\Ga$ is a cusp; compare 1.3, (A).}

(II) Assume that $\Oh_C=\prod R_E/(t_E{}^2)$, with all the $n_E=2$, and
that $R_E/(t_E)=k(E)=K$ for each $E$; write $T^*_E=(t_E)/(t_E{}^2)\subset
R_E/(t_E{}^2)$ for the cotangent space of\/ $R_E$, a $1$@-dimensional
$K$@-vector space. Let $\Oh_D\subset\Oh_C$ be a half-filling. Then
$m_D\subset\sum T^*_E$ is a $K$@-vector subspace of codimension $1$, and if
there are $\ge2$ components $E$, it involves every summand:
 $$
m_D\not\subset\sum_{E'\ne E} T^*_{E'}.
 $$

That is (for varieties over a field, for simplicity), suppose that $X$
is nonnormal and Gorenstein, and that every component $E$ has multiplicity $2$
in $C$, and maps birationally to $\Ga=D\red\subset\Sing X$. Then the
transverse singularity of $X$ along $\Ga$ is a cusp $(y^2=x^3)$ if $r=1$, a
tac\-node $(y^2=x^2y)$ if $r=2$, and $r$ concurrent lines in $\aff^{(r-1)}$
with
no $(r-1)$ in a hyperplane for $r\ge 3$.

{\rm In (II), I assume that the residue field extensions are trivial
$R_E/t_E=k(E)=K$ mainly out of spinelessness; this case is sufficient for my
del
Pezzo surfaces, since by Theorem~1.1, the only multiple locuses in $C$ are
double lines. More generally, $T_E=(t_E)/(t_E{}^2)\cong k(E)$ and
$m_D=\ker\psi$, where $\psi\:\sum T_E\to K$ is a nonzero multiple of the trace
map $\Tr_{k(E)/K}$ on each piece.}\endproclaim

\demo{{\rm3.8} Proof} (I) The notation $\Oh_C=\prod A_E$ with
$A_E=R_E/(t_E{}^{n_E})$ is as above; since $X$ is nonnormal, $\sum
n_E\cdot[k(E):K]\ge2$, so that either there's more than one component $E$, or a
component with $n_E\ge2$ or $[k(E):K]\ge2$.

The dual statement to the local condition $\dim\Hom_{\Oh_D}(\Oh_D,K)=1$ on
$\Oh_D$ (see 3.5, (b)) is that its {\it socle}
 $$
\{s\in\w_D\mid m_D\cdot s=0\}=\Hom_{\Oh_D}(K,\w_D)
 $$
is $1$@-dimensional over $K$. Assuming $\Oh_D\subset\Oh_C$ is a
half-filling, $\w_D=\Oh_C/\Oh_D$.

If a component $E$ of $C$ has $n_E=1$ then $A_E=k(E)\subset\Oh_C$ is killed by
$m_D$; hence $S=\sum_{n_E=1}k(E)\subset\Oh_C $ is a $K$@-vector subspace of the
socle of $\Oh_C$ (as an $\Oh_D$@-module). Consider the projection of $S$ to
$\Oh_C/\Oh_D$. The local condition means that
 $$
\Oh_D\cap S\cong
\cases
K&\text{if there are no further components;}\\
0&\text{if there are components with $n_E\ge2$.}
\endcases
 $$
Therefore, if there are no further components, $S$ contributes
 $$
\sum_{n_E=1}[k(E):K]-1\le1
 $$
to the $K$@-dimension of the socle, which leads at
once to the conclusion (I).

If there are components with $n_E\ge2$, I look for a contradiction. First, $S$
contributes $\sum_{n_E=1}[k(E):K]\le1$ to the dimension of the socle. It
follows
that $S=K$, hence $\Oh_C$ is of the form $K\times\Oh_{C'}$, and by the local
condition $\Oh_D$ doesn't meet the first factor, so projects isomorphically to
$\Oh_D\subset\Oh_{C'}$. But $\ell(\Oh_D)>\ell(\Oh_{C'}/\Oh_D)$, so
$\Oh_D\subset\Oh_{C'}$ is not a part-filling by 3.6, Step~4; thus
$\Oh_{C'}/\Oh_D$ is an unfaithful $\Oh_D$@-module. If $n\in\Oh_D$ kills
$\Oh_{C'}/\Oh_D$ then there is a case division:

\smallskip\noindent
{\bf Subcase} $n$ a unit. Then $\Oh_{C'}/\Oh_D=0$, which gives $\Oh_C\cong
K\times K$; this certainly happens, but not under the current case assumption
$n_E\ge2$.

\smallskip\noindent
{\bf Subcase} $n\in m_D$. Then $n\cdot\Oh_{C'}\subset\Oh_D$ and $n\cdot K=0$
implies at once that
$n\cdot\Oh_C=n\cdot(K\times\Oh_{C'})\subset\Oh_D$, contradicting the
assumption that $\Oh_C/\Oh_D$ is faithful.

(II) Suppose there are $r$ factors so that $\Oh_C=\prod R_E/(t_E{}^2)$ has
length $2r$. Thus $\Oh_D$ has length $r$, and by the local condition,
$\Oh_D\subset\Oh_C$ maps to the diagonal $K\subset\prod k(E)$ with kernel
$m_D=\Oh_D\cap\sum T_E$ of length $(r-1)$. If $r=1$ there's nothing more to
prove.

The dual of the diagonal map of $\Oh_D\subset\Oh_C$ to $K\subset\prod k(E)$ is
the inclusion of the socle of $\w_D=\Oh_C/\Oh_D$, which is a map $K\to(\sum
T_E)/m_D$ with each component nonzero; by duality this is surjective, so that
$m_D\subset\sum T_E$ is a codimension one subspace involving each summand.

The geometric statement about $Y\to X$ is easy to see. For example if $r=1$
then $C=2E$, there is a morphism $\fie\:2E\to\Ga$ which is birational when
restricted to $E$, and $\Oh_X\subset\Oh_Y$ is the ring of functions $f$ such
that $f\rest{2E}\in\Oh_\Ga$. That is, $Y$ is pinched along $\fie$ to give a
codimension 1 locus of cusps. \QED\enddemo

\subheading{{\rm3.9}\enspace Derivations and case (II) of Theorem~3.7} I now
show how a subring $\Oh_D\subset\Oh_C$ in case (II) of the Theorem~3.7 is
specified in practical calculations. Let $\pi\:Y\to X$ be a normalisation of a
reduced $S_2$ scheme with $\fie\:C\to D$ the conductor locus. Suppose that
$\Ga=D\red$ is irreducible, and that above its generic point, $\pi\:Y\to X$
falls under case (II) of Theorem~3.7.

The glueing map $\fie\:C\to D$ can be factored as a composite of two maps; let
$D^+$ be the variety homeomorphic to $D$, but with
$\Oh_{D^+}\subset\fie_*\Oh_C$ defined by
 $$
\Oh_{D^+}=\bigl\{\{h_i\}\in\fie_*\Oh_C\bigm|\im h_i\in k(E)=K\text{ is
independent of $i$}\bigr\}.
 $$
That is, $C\to D^+$ does only the set-theoretic glueing, leaving the $r$
tangent spaces generically transversal, whereas $D^+\to D$ squashes up a first
order disc in $r$@-space over $k(\Ga)$ into a first order disc in
$(r-1)$@-space. Geometrically, $D^+$ corresponds to the ``transversalisation''
of the nonnormal locus of $X$, that is, $Y\to X^+\to X$ where $X^+$ has a
multiple curve with $r$ generically transverse branches.

Although this is not very intuitive, it's important to understand that already
on the generic point, $D^+\to D$ contains nontrivial information of a
differential nature: for example, if $r=1$ then $C=D^+$, and $\fie$ maps a
normal field to $C\red\subset C$ to a vector field on $D$.

\proclaim{Proposition} Under the above assumptions, there exists a
(rational) derivation
 $$
\De\:\Oh_{D^+}\to(\Oh_{D^+\red})\gen=k(E)
 $$
such that $\Oh_D=\ker\De\subset\Oh_{D^+}$.\endproclaim

See 4.4--5 for a much more concrete description in the particular case of
interest.

\demo{Proof} Write $N_1$ and $N_2$ for the nilpotent kernels of
$\Oh_{D^+}\to\Oh_{D^+\red}$ and $\Oh_{D}\to\Oh_{D\red}$; this gives the exact
diagram
 $$
\spreadmatrixlines{4pt}
\matrix\format &\,\c\,\\
0&\to&N_2&\to&\Oh_D&\to&\Oh_{D\red}&\to& 0\\
&&\bigcap&&\bigcap&&\bigcap\\
0&\to&N_1&\to&\Oh_{D^+}&\to&\Oh_{D^+\red}&\to& 0\\ \vspace{-.8ex}
&&\downarrow&&\downarrow\\
0&\to&N_1/N_2&\to&\Oh_{D^+}/\Oh_D.\\
\endmatrix
 $$

Looking only at generic stalks, the assumptions in Theorem~3.7, (II) give that
$\Oh_{D\red}=\Oh_{D^+\red}$, so that
 $$
(N_1/N_2)\gen=(\Oh_{D^+}/\Oh_{D})\gen,
 $$
and the left-hand side is a 1@-dimensional vector space over $k(E)$.

Now (a) all the above maps are $\Oh_D$ linear; and (b) since
$N_1\subset\Oh_{D^+}$ is an ideal with square zero, all the sheaves $N_1$,
$N_2$ and $N_1/N_2$ in the left-hand column have both $\Oh_{D^+}$@-module and
$\Oh_{D}$@-module structures that are compatible. The fact that
$\De\:\Oh_{D^+}\to(\Oh_{D\red})\gen$ is well defined and a derivation
follows formally. \QED\enddemo

\heading\S4. The glueing map $\fie\:C\to D$ and proof of Theorem~1.5
\endheading

\myno{4.0} The starting point in this section is a disjoint union
 $$
(C\subset Y)=\coprod (C_i\subset Y_i)
 $$
of surfaces $C_i\subset Y_i$ taken from the list of Theorem~1.1 (all defined
over the same algebraically closed field). I study the possibilities for
glueing
$Y$ by a morphism $\fie\:C\to D$ to get a connected Gorenstein surface $X$;
since
$\pi^*\w_X=\w_Y(C)=\Oh_Y(-1)$, it will follow automatically that the invertible
sheaf $\w_X{}^{-1}$ is ample, so $X$ is a del Pezzo surface. The main result
Theorem~1.5 can be viewed as a classification of nonnormal del Pezzo surfaces:
under extra conditions they are just the surfaces $X$ described in 1.3--4.

A nice simplifying feature is that I hardly need to work with $X$ and $Y$ at
all: given that $\w_Y(C)\cong\Oh_Y(-1)$, it follows by Corollary~2.8, iv that
$X$ is Gorenstein if and only if $\ker\{\Tr_{C/D}\:\fie_*\w_C\to\w_D\}$ is an
invertible $\Oh_D$@-module. Technically, the main point is to determine the map
$\Tr_{C/D}$, in order to verify the conditions on $s\in\fie_*\w_C$ that it (a)
belongs to $\ker\Tr$, and (b) is an $\Oh_D$@-basis of $\ker\Tr$.

\myno{4.1} The first approximation to Theorem~1.5 is a numerical treatment that
considers only the degrees of the curve components of $\Sing X$ and $C$ with
respect to the polarisations $\Oh_X(1)$, $\Oh_Y(1)$, and the nature of $X$ at
the generic point. Note that $\deg\Ga=1$ or $2$ doesn't imply automatically
that $\Ga$ is isomorphic to a line or conic, since a priori $\Oh_X(1)$ is not
very ample.

\proclaim{Lemma} One of the following holds (compare 1.3--4):

\myno{(A)} $\Sing X=\Ga$ is irreducible with $\deg\Ga=1$, and
$C=\pi^{-1}\Ga\subset Y$ a nonsingular conic with
$\deg(\fie\:C\to\Ga)=2$; in this case $Y$ is irreducible, and $X$ has
ordinary double points in codimension $1$ along $\Ga$. This includes
the possibility of inseparable ordinary double points in characteristic $2$,
see
the note in Theorem~3.7, (I).

\myno{(B)} $\Sing X=\Ga$ is irreducible with $\deg\Ga=2$; then $Y$
has exactly $2$ components $Y_1$, $Y_2$, the curves $C_1\subset Y_1$ and
$C_2\subset Y_2$ are nonsingular conics mapping birationally to $\Ga$, and
$X$ has ordinary double points in codimension $1$ along $\Ga$.

\myno{(C$_1$)} $\Sing X=\Ga$ is irreducible with $\deg\Ga=1$ and
$C=\pi^{-1}\Ga\subset Y$ a line pair; in this case $Y$ is irreducible,
and $X$ has ordinary double points in codimension $1$ along $\Ga$.

\myno{(C$_r$)} $X=\bigcup_{i=1}^r X_i$ is a cycle of $r\ge2$ components
$X_i$, where $X_i$ and $X_{i+1}$ meet generically transversally along a curve
$\Ga_i$ with $\deg\Ga_i=1$ (the indexes are taken cyclically, that
is, $X_{r+1}=X_1$). Each conductor locus $C_i\subset Y_i$ is a line pair
$C_i=\ell_i\cup\ell'_i$, whose two components map birationally to $\Ga_{i-1}$
and $\Ga_i$.

\myno{(D)} $\Sing X=\Ga$ is irreducible with $\deg\Ga=1$, and the
conductor locus $C_i\subset Y_i$ in each component of\/ $Y$ is a double
line; in this case $X$ has the singularities in codimension $1$ along
$\Ga$ described in Theorem~3.7, (II).\endproclaim

\demo{Proof} Each $C_i$ is a conic, so is either reduced or a double line.
Above a component $\Ga\subset\Sing X$, the behaviour of $Y\to X$ is described
by Theorem~3.7, (I) if $\pi^{-1}\Ga\subset C$ is reduced, or Theorem~3.7, (II)
if it consists of double lines.

Because the glueing takes place in codimension $1$ by Proposition~2.2, in
order for $X$ to be connected, the components of $Y$ must be joined together
along components of $C_i$. If $C$ has a nonsingular conic component, then
$C\to\Ga$ is either a double cover of a curve of degree $1$, which implies
that $Y$ is connected; or two conics covering a curve of degree $2$, which
implies that $Y$ has exactly $2$ components. This gives cases (A--B). If say
$C_1=\ell_1\cup\ell'_1\subset Y_1$ is a line pair, then $\fie\:C\to D$ must
glue $\ell'_1$ birationally to some other line, say $\ell''_1=\ell_2\subset
C_2\subset Y_2$, and so on, to form a cycle, giving case (C$_r$).

The remaining possibility is that $C\to D$ is made up of $r$ double lines
mapping
generically to $D\subset X$ as in Theorem~3.7, (II). \QED\enddemo

\subheading{{\rm4.2}\enspace Proof of Theorem~1.5 for reduced $C$} Case (A) is
very easy: the conductor $\Cee\subset X$ defines the reduced curve $\Ga$ at its
generic point, hence everywhere. Therefore $D=\Ga$. Since $C$ is normal, the
form ``$D$ normal under $C$'' of the $S_2$ condition discussed in
Proposition~2.2 implies that $D$ is normal, so $D\cong\proj^1$, and
$\fie\:C\to D$ is isomorphic to a linear projection of a conic to a line.

Case (B) is similar: each of the two components $C_1$, $C_2$ maps birationally
to $D=\Ga$ under $\fie$; so $S_2$ again implies that $D$ is normal, and the
glueing just consists of identifying $C_1$ and $C_2$ by an isomorphism to $D$.

Case (C): as before, $\Cee\subset\Oh_X$ defines the reduced curve
$D=\bigcup\Ga_i$ at each generic point, hence everywhere, so that $D$ is
reduced. According to Lemma~4.1, every component $\Ga_i$ has exactly two
line components $\ell'_i$ and $\ell''_i$ of $C$ mapping birationally to it. At
last I can use the Gorenstein condition of Corollary~2.8, iii to prove
something nontrivial:

\proclaim{{\rm4.3}\enspace Claim} There is a unique vertex point $P\in D$ such
that
 $$
\fie^{-1}P=\{\text{\rm nodes of $C$}\},
 $$
in other words, all the nodes of the $C_i$, and no other points, map to $P$.
\endproclaim

The point is the following: $\fie$ identifies two lines $\ell'_i$ and
$\ell''_i=\ell_{i+1}$ birationally to $\Ga_i$, and each of these has a marked
point $P'_i\in\ell'_i$ and $P''_i\in\ell''_i$, namely the nodes of
$C_i=\ell_i\cup\ell'_i$ and $C_{i+1}=\ell_{i+1}\cup\ell'_{i+1}$. The claim is
that these two points match up under the identification. Assuming this, it's
clear from the $S_2$ condition that the $r$ components $\Ga_i$ are all
normal, and that they define $r$ linearly independent directions in the
tangent space to $D$ at $P$. This implies Theorem~1.5 in case (C).

\demo{Proof of Claim~4.3} As pointed out in Remark~2.9, the trace map is
birational in nature. The birational identifications of $\ell'_i$ and
$\ell''_i$ with $\Ga_i$ under $\fie$ identifies the 3 generic stalks of the
dualising sheaves, $\w_{\ell',}{}\gen$, $\w_{\ell'',}{}\gen$ and
$\w_{\Ga,}{}\gen$; here the subscript $\gen$ denotes the generic stalk, and
I omit the secondary subscript $i$ from now on. The calculation of $\Tr_{C/D}$
is trivial, namely
 $$
\w_{\ell',}{}\gen\oplus\w_{\ell'',}{}\gen\ni(s',s'')\mapsto
s'+s''\in\w_{\Ga,}{}\gen.
 $$
Thus $(s',s'')\in\ker\Tr$ if and only if $s'=-s''$.

The birational identification of the nonsingular curves $\ell'_i$ and
$\ell''_i$ extends to a biregular identification. The elements
$s'\in\w_{\ell',}{}\gen$ and $s''\in\w_{\ell'',}{}\gen$ are rational sections
of $\w_{\ell'}$ and $\w_{\ell''}$, and if $s'=-s''$ then obviously the zeros
and poles of $s'$ and $s''$ must match up under the identification. To get a
contradiction, suppose that $Q'\in\ell'_i$ and $Q''\in\ell''_i$ are
identified, where $Q'=\ell_i\cap\ell'_i\in C_i$ is the node and $Q''\in
C_{i+1}$ a nonsingular point. The Gorenstein condition of Corollary~2.8, iii
is that $\ker\Tr$ contains an $\Oh_C$@-basis of $\w_C$ at each point of $C$. I
will deduce from this that $\ker\Tr$ has a section $(s',s'')$ over the generic
point of $\Ga_i$ such that $s''$ is regular at $Q''$, but $s'$ has a pole at
$Q'$, contradicting what I just said.

Let $(s',s'')\in\w_{\ell',}{}\gen\oplus\w_{\ell'',}{}\gen$ be a section of
$\fie_*\w_C$ at the generic point of $\Ga_i$. First, $Q''$ is not a
node, so $\ell''_i=C$ near $Q''$, and in order for $(s',s'')$ to be in
$\fie_*\w_C$ near $\fie(Q'')$, clearly $s''$ must be regular at $Q''$.
But on the other hand, if $Q'=\ell_i\cap\ell'_i\in C_i$ is the node of $C_i$,
then an $\Oh_{C_i}$@-basis of the stalk of $\w_{C_i}$ at $Q'$ is given by
$(s,s')$ where $s\in\w_{\ell_i}(Q')$ and $s'\in\w_{\ell'}(Q')$ are
bases; that is, $s'$ must have a pole at $Q'$. This completes the
contradiction. \QED\enddemo

\subheading{{\rm4.4}\enspace The nonreduced case} I first describe without
proof the argument in the most important case $r=1$ so that the trusting or
exhausted reader can skip the rest of the paper. By Lemma~4.1, in this case
$C$ is a double line $C\cong2\ell\subset\proj^2$, and $\fie\:C\to D=\Ga$ a
morphism whose restriction to the line $\ell=C\red$ is birational. The $S_2$
condition that $\fie_*\Oh_C/\Oh_\Ga$ is torsion-free does not imply $\Ga$
nonsingular.

An affine piece of $C$ is $(y^2=0)\subset\aff^2_{x,y}$, so that the generic
stalk of $\Oh_C$ is $k(x)[y]/(y^2)$, and $k(\Ga)=k(x)$; the map $\fie\:C\to
\Ga$ at the generic point corresponds to an inclusion
 $$
\spreadmatrixlines{3pt}
\matrix\format\r&\c&\l\\
k(x)&{}\hookrightarrow{}&k(x)[y]/(y^2)\\
\text{by}\quad f&\mapsto&f-h_\fie f'y\\
\endmatrix,
\qquad\text{where}\;f'={{\dd f}\over{\dd x}}\;\text{and}\;
h_\fie=h_\fie(x)\in k(x);
 $$
$h=h_\fie$ can be arbitrary. In other words, generically $\fie$ is a
``projection'' of the double line $C=2\ell$ back to the reduced line
$\ell\bir\Ga$; and the projection of a ``normal field'' to $\ell\subset C$
defines a derivation $f\mapsto hf'$ of $\Oh_\ell$, that is, a rational vector
field on $\ell$.

In this situation, it's clear from the $S_2$ condition that $\ell\to\Ga$ is
one-to-one, and that $\Oh_\Ga\subset\Oh_C$ is the sheaf of all
functions of the form $f-hf'y\in\Oh_C$, where both $f$ and $hf'$ are regular on
$\ell$; that is,
 $$
\Oh_\Ga\cong\{f\in\Oh_\ell\mid hf'\in\Oh_\ell\};
 $$
thus $\Oh_\Ga\subsetneq\Oh_\ell$ at the poles of $h$, so $\Ga$ is nonnormal
there. Clearly $\xi=x-hy$ is a rational section of $\Oh_\Ga$ that bases the
function field, $k(\Ga)=k(\xi)$. One calculates (see 4.8) that in terms of the
natural bases
 $$
s_C=\Res_{\aff^2/C} \left({{\dd x\wedge\dd
y}\over{y^2}}\right)\in\w_C\quad\text{and}\quad\dd \xi\in\w_{k(\Ga)},
 $$
the trace map $\Tr_{C/D}\:\w_C\to\w_D$ at the generic point of $C$ is
 $$
\Tr_{C/D}\bigl((f+gy)\cdot s_C\bigr)= \bigl(g+(hf)'\bigr)\cdot\dd\xi,
 $$
so that $(f+gy)\cdot s_C\in\ker\Tr\iff g=-(hf)'$.

For $(f+gy)\cdot s_C$ to be a basis of $\w_C$ at a point $P\in C$ it's
necessary and sufficient for $f$, $g$ to be regular and $f(P)\ne0$. Now it's
easy to see that at $P\in\ell$, there exists a unit $f\in\Oh_\ell$ with
$(hf)'\in\Oh_\ell$ if and only if either $h$ is regular at $P$, or $\cha k=p$
and $h$ has a pole at $P$ of order divisible by $p$ (compare the proof of
Theorem~4.6). Using the Gorenstein criterion Corollary~2.8, iii, the results
stated in the introduction follow from this: if $\cha k=0$ then $C\to D$
must be a linear projection of a conic to a line. But in characteristic $p$,
the rational function $h$ is allowed to have poles of order $=np$ for any
$n\ge1$; at such a pole $P\in\ell$, it's clear that the local ring is of the
form
 $$
\Oh_{\Ga,P}\cong k[x^i\mid\{i\equiv 0\text{ mod $p$, or }i\ge np\}]_{(0)}
 $$
(localised at $x=0$), so that $\Ga=D$ has a ``wild'' cusp. Arbitrary poles of
this form can happen, so $H^1(\Oh_X)$ can be arbitrarily large. Clearly the
local ring of $\Ga$ at $P$ needs $p$ generators $x^p$ and
$x^{np+1},\dots,x^{np+p-1}$, so that this type of singularity cannot occur if
$\Ga$ is contained in a smooth $3$@-fold and $p\ge5$ (or, more generally, if
it's given that $\dim T_PX\le\cha k$). In fact it can be checked (see
Exercise~4.12) that in the case when $Y$ is nonsingular,
 $$
\cases
 \dim T_PX=\dim T_P\Ga=p&\text{if $p\ge3$;}\\
 \dim T_P\Ga=2,\dim T_PX=3&\text{if $p=2$.}\\
\endcases
 $$
Thus in this case, $X$ has hypersurface singularities if $\cha k=2$ or $3$.

\myno{4.5} In the general case, $C=\coprod_{i=1}^rC_i$ with each $C_i$ a double
line $C_i=2\ell_i\subset\proj^2$. Restricting $\fie$ defines a birational map
$\fie_i\:\ell_i\to\Ga=D\red=\Sing X$ on each of the reduced curves; I write
$\ell$ for the normalisation of $\Ga$, so that the $\fie_i$ define
isomorphisms $\ell_i\cong\ell$. Now, as in 3.9, let $\fie^+\:C\to D^+$ be the
morphism that glues the $C_i$ to each other along this isomorphism. In other
words, $D^+$ is isomorphic to a first order infinitesimal neighbourhood of
$\ell=\proj^1\subset\proj^{r+1}$, and the $C_i\subset D^+$ are double lines in
$r$ transverse planes through $\ell$. Let $x$ be an affine parameter on $\ell$
and $y_i$ a linear form in the plane of $C_i$ vanishing along $\ell_i$.

 $$\spreadmatrixlines{4pt}\matrix\format\r&\c\qquad&\l\\
&\Oh_C&C=\coprod C_i\\
&\bigcup&\kern2pt\downarrow\fie^+\\
\Oh_\ell\twoheadleftarrow&\Oh_{D^+}&D^+\supset\ell\\
&\bigcup&\kern2pt\downarrow\\
\Oh_\Ga\twoheadleftarrow&\Oh_D&D\supset\Ga.\\
\endmatrix
 $$

Obviously, the generic stalk of $\Oh_{D^+}$ is $k(x)[\underline y]/(\underline
y)^2$,
where $\underline y$ is short for $y_1,\dots,y_r$. By Proposition~3.9,
$\Oh_D=\ker\De$ where $\De\:\Oh_{D^+}\to k(x)$ is a rational derivation. But
any derivation $\De\:k(x)[\underline y]/(\underline y)^2\to k(x)$ over $k$ is
of the form
 $$
\De(a,\underline b)\:f(x)+\sum g_i(x)y_i\mapsto af'+\sum b_ig_i
\qquad\text{where }\ f'={\dd f\over\dd x},
 $$
for some $a, b_1,\dots,b_r\in k(x)$. Theorem~3.7, (II) contains the assertion
that $b_i\ne0$ for each $i$; the case $a=0$ is most welcome. A change of basis
in $k(x)$ multiplies through $(a,\underline b)$ by a nonzero element of $k(x)$.

\myno{4.6} I choose coordinates $x_i$, $y_i$ on an affine piece of each $C_i$,
where $x_i$, $y_i$ are the pullback of the coordinates $x$, $y_i$ on $D^+$; a
function on $C_i$ is of the form $f_i+g_iy_i$ where $f_i$, $g_i$ are rational
functions of $x_i$, and from now on I usually omit mention of the substitution
$x_i\mapsto x$ corresponding to the (nameless) identification of $\ell_i$ and
$\ell$.

\proclaim{Theorem} $\ker\Tr_{C/D}$ contains a basis of $\fie_*\w_C$ at $P\in
\ell$ (compare Corollary~2.8, iii) if and only if\/ $b_i/b_j$ is a unit of\/
$\Oh_{\ell,P}$ for all\/ $i$, $j$, and either each $a/b_i$ is regular at $P$,
or $\cha k=p$ and $a/b_i$ has a pole at $P$ of order divisible by
$p$.\endproclaim

\demo{Proof} An element of $\fie_*\w_C$ is of the form
 $$
s(\underline f,\underline g)=\sum_{i=1}^r (f_i+g_iy_i)\cdot s_i,
\quad\text{where }\
s_i=\Res_{\aff^2/C_i}{{\dd x_i\wedge\dd y_i}\over{y_i{}^2}}
 $$
is a basis of $\w_{C_i}$ over $\Oh_{C_i}$, and $\underline f=(f_1,\dots,f_r)$,
$\underline g=(g_1,\dots,g_r)$. The calculation of $\Tr_{C/D}$ given in 4.8
below proves that
 $$
s(\underline f,\underline g)\in\ker\Tr_{C/D}\enspace\iff\enspace{f_i\over
b_i}={f_1\over b_1}\enspace\text{for all $i\ge2$,\enspace and }\left({af_1\over
b_1}\right)'=-\sum_{i=1}^r g_i.\tag1
 $$

Assuming (1) the result is easy: in order for $s(\underline f,\underline g)$
to base $\w_{C_i}$ at each $P_i$ lying over $P\in\ell$, each $f_i$ must be a
unit, and $g_i$ regular. Therefore $b_i/b_j=f_i/f_j$ is a unit for all $i$,
$j$. And I must be able to solve
 $$
\left({af_1\over b_1}\right)'=-\sum g_i\qquad\text{with $g_i$ regular and $f_1$
a unit}.\tag2
 $$
Suppose that $x$ is a local parameter at $P\in\ell$, and that $a/b_1=x^m\cdot
u$ with $u$ a unit. It's obvious that regular functions $g_i$ can be chosen to
solve (2) if and only if the left-hand side is regular. The left-hand side
of (2) is the derivative of $f_1\cdot u\cdot x^m$, which up to multiplying
by a nonzero constant is $mx^{m-1}+\text{higher order terms}$. This is not
regular if $m<0$ and $m\ne0\in k$, so that the condition on $a/b_1$ is
necessary. When $m<0$, I can choose $f_1$ such that $f_1$ is a unit and
$f_1\cdot u-1$ has a zero of order $\ge-m$ at $P$, so that
 $$
\left({af_1\over b_1}\right)'=mx^{m-1}+\text{regular terms};
 $$
if also $m=0\in k$ then this is regular, so I can solve (2) as before. This
proves the theorem, assuming (1).\enddemo

\subheading{{\rm4.7}\enspace Tame} Using the notation of Theorem~4.6, $D\red$
is nonsingular if and only if $D\red\cong\ell=D^+\red$ under each $P\in\ell$,
which happens if and only if $a/b_i$ is regular, so that
$x-(a/b_i)y_i\in\Oh_D$ is a local parameter. In local terms (also globally,
see 4.11), if $\cha k=p$ then $a/b_i$ can have poles of order $=np$ an
arbitrary multiple of $p$; as discussed at the end of 4.4, it's easy to see
that in this case the local ring of $D\red$ is $k[x^i\mid\{i\equiv 0\text{ mod
$p$, or }i\ge np\}]_P$, so that it requires $p$ generators.

With $\pi\:Y\to X$ and $\fie\:C\to D$ as in (0.2), I say that
$\fie\:C\to D$ (or $X$ itself) is {\it tame} if either $C$ is reduced, or $C$
is nonreduced and $D\red$ nonsingular. Tameness is automatic under any of the
following alternative extra assumptions on $X$:

(i) The double locus $C\subset Y$ has a reduced component (see Theorem~3.7,
(I)).

(ii) $\cha k=0$.

(iii) $\cha k=p\ge5$ and $X$ is locally a divisor in a nonsingular 3@-fold,
or, more generally, $\cha k=p>\dim T_PX$ for all $P\in X$.

(iv) $C=\coprod C_i$ with each $C_i$ a double line and $\Sing
X=D\red\cong\proj^1$ (by definition).

(v) $\Oh_X(1)\rest{D\red}$ is very ample (since then every component of
${D\red}$ is isomorphic to $\proj^1$ or a plane conic).

(vi) $H^1(\Oh_X)=0$ or $\chi(\Oh_X)=1$ (see Corollary~4.10).

\subheading{{\rm4.8}\enspace The calculation of $\Tr_{C/D}\:\fie_*\w_C\to\w_D$}
The
unpleasant thing is that $\Tr_{C/D}$ is $\Oh_D$@-linear, but not $k(x)$@-linear
in
general; to calculate it, I introduce new coordinates at the generic points of
$C$ that are rational sections of $\Oh_D\subset\fie_*\Oh_C$. The whole point is
to write down sections of $\w_C$ in terms of these new coordinates.

Recall that by Proposition~3.9 and 4.5, $\Oh_D\subset\Oh_{D^+}$ is the
subring defined by
 $$
\Oh_D=\bigl\{f+\sum g_iy_i\in\Oh_{D^+}\bigm| af'+\sum b_ig_i=0\bigr\}.
 $$
So in particular
 $$
\split
\xi&=x-(a/b_1)y_1\\
\text{and}\quad\eta_i&=y_i-(b_i/b_1)y_1\quad\text{for $i=2,\dots,r$}
\endsplit
 $$
are rational sections of $\Oh_D$. I've chosen $\xi$ and $\eta_i$ so that they
only mess up $C_1$, and not $C_i$ for $i=2,\dots,r$: the pullback of $\xi$ to
$C_1$ is $\xi_1=x_1-(a/b_1)y_1$, and to $C_i$ is $\xi_i=x_i$ for $i\ge2$;
whereas each $\eta_i$ for $i\ge2$ pulls back to $y_i$ on $C_i$, to
$-(b_i/b_1)y_1$ on $C_1$ and to $0$ on all components $C_j$ with $j\ne i$.

I have to translate
 $$
s(\underline f,\underline g)=\sum_{i=1}^r\bigl(f_i(x_i)+g_i(x_i)y_i\bigr)\cdot
s_i\in\fie_*\w_C
 $$
into the new coordinates. For this, first define new generic bases
 $$
s_i'=\Res_{\aff^2/C_1}{{\dd\xi_i\wedge\dd y_i}\over{y_i{}^2}}\in\w_{C_i}.
 $$
There's not very much to calculate for the terms with $i\ge2$, since $s_i=s'_i$
and $f_i(x_i)+g_i(x_i)y_i=f_i(\xi_i)+g_i(\xi_i)y_i$. The fun comes when $i=1$;
note that $\dd x_1=\bigl(1+(a/b_1)'y_1\bigr)\cdot\dd\xi_1+(a/b_1)\dd y_1
\vphantom{\Big(}$. Thus
 $$
s_1=\Res_{\aff^2/C_1}{{\dd x_1\wedge\dd y_1}\over{y_1{}^2}}=
\bigl(1+(a/b_1)'y_1\bigr)\cdot s_1'.
 $$
Now using $y_1{}^2=0$ and the Taylor expansion for
$f_1(x_1)=f_1(\xi_1+(a/b_1)y_1)$ gives
 $$
\split
&\bigl(f_1(x_1)+g_1(x_1)y_1\bigr)\cdot s_1\\
&\qquad=\bigl(f_1(\xi_1)+(a/b_1)f'y_1+g_1y_1\bigr)\cdot
\bigl(1+(a/b_1)'y_1\bigr)\cdot s_1'\\
&\qquad=\bigl(f_1(\xi_1)+(af/b_1)'y_1+g_1y_1\bigr)\cdot s_1'.\\
\endsplit
 $$
Therefore,
 $$
s(\underline f,\underline g)=\left({af\over b_1}\right)'y_1\cdot s_1'
+\sum_{i=1}^r\bigl(f_i(\xi_i)+g_iy_i\bigr)\cdot s_i'.\tag3
 $$

\demo{Proof of $(1)$} Since $\xi$ is a rational section of $\Oh_D$, the trace
$\Tr_{C/D}$ on generic sections is $k(\xi)$@-linear. Hence I view the generic
stalk of $\fie_*\w_C$ as $\Hom_{k(\xi)}(\fie_*\Oh_C, k(\xi))$. The trace is
just the restriction to $\Oh_D$ of elements of $\Hom_{k(\xi)}(\fie_*\Oh_C,
k(\xi))$, so a generic section of $\fie_*\w_C$ is in the kernel of $\Tr_{C/D}$
if and only if it kills the elements $1\in\Oh_D$ and $\eta_i\in\Oh_D$ for
$i\ge2$, that form a $k(\xi)$@-basis of $(\Oh_D)\gen$.

Now making the usual identification of $\w_{C_i}$ with differentials, I get
that the element $\bigl(f_i(\xi_i)+g_i(\xi_i)y_i\bigr)\cdot s_i'$ is the
function on $\Oh_{C_i}$ taking $a(\xi_i)+b(\xi_i)y_i\in\Oh_{C_i}$ into the
residue of the product, that is, into the coefficient $ag_i+bf_i\in k(\xi)$ of
$y_i$. Now using (3), it's easy to see that evaluating $s(\underline
f,\underline g)$ on $1_{\Oh_D}=\{1,\dots,1\}\in\fie_*\Oh_C$ and
$\eta_i=y_i-(b_i/b_1)y_1$ gives
 $$
\left({af\over b_1}\right)'+\sum_{i=1}^r g_i\qquad\text{and}\qquad
f_i-(b_i/b_1)f_1.
 $$
The conditions for $\Tr_{C/D}s(\underline f,\underline g)=0$ are obtained by
setting both of these to $0$. This proves (1) and completes the proof of
Theorem~4.6. \QED\enddemo

\proclaim{{\rm4.9}\enspace Proposition} Write $\ell,\Oh_\ell(1)$ for $\proj^1$.
Then
in the notation of 4.5, the structure sheaf of\/ $D^+$ is isomorphic
to $\Oh_{D^+}\cong\Oh_\ell\oplus\bigoplus_{i=1}^r\Oh_\ell(-1)$ as a sheaf of
rings, with the second summand $N_1$ an ideal of square zero.

Set $N_2=\ker\{\Oh_D\to\Oh_{D\red}\}$; then $N_2$ has an $\Oh_\ell$@-module
structure, and, assuming that $X$ is Gorenstein,
$N_2\cong\bigoplus_{i=1}^{r-1}\Oh_\ell(-1)$. Suppose in addition that $X$ is
tame (see 4.7). Then as a sheaf of rings
$\Oh_D\cong\Oh_\ell\oplus\bigoplus_{i=1}^{r-1}\Oh_\ell(-1)$. In other words,
$D,\Oh_D(1)$ is isomorphic to a first order neighbourhood of\/ $\proj^1$ in
$\proj^r$, with the morphism $\fie\:\coprod C_i\to D$ linear on each
component.\endproclaim

\demo{Proof} The first part is easy, since $\Oh_{D^+}$ fits in the locally
split exact sequence $0\to N_1\to\Oh_{D^+}\to\Oh_\ell\to0,$ with
$N_1=\bigoplus_{i=1}^r\Oh_\ell(-1)$, any two local splittings $s_i$ and $s_j$
differ by a derivation $\de_{i,j}=s_i-s_j\:\Oh_\ell\to N_1$, and
$H^1(\CDer(\Oh_\ell,N_1))=H^1(\CHom(\Omega^1_\ell,N_1))=0$.

The $\Oh_{D^+}$ or $\Oh_\ell$@-module structure of $N_2$ comes from the fact
that $N_1$ is nilpotent of square $0$. Consider the projection $N_1\to
N_1/N_2$. By the $S_2$ condition $N_1/N_2$ is a line bundle over $\ell$, and
$N_1\to N_1/N_2$ is the map defined by the $b_i$. Since by Theorem~4.6 the
$b_i/b_j$ are units, it follows that each direct summand $\Oh_\ell(-1)$ of
$N_1$ (corresponding to the $y_i$) maps isomorphically to $N_1/N_2$, so that
$N_1/N_2\cong\Oh_\ell(-1)$. Thus the kernel $N_2$ is also a direct sum of
copies of $\Oh_\ell(-1)$.

In the notation of Theorem~4.6, if any of the $a/b_i$ is a unit, then I can
choose $\xi=x-b_iy_i\in\Oh_D\subset\Oh_{D^+}$ which maps to a local parameter
of $\Oh_\ell$; since I'm assuming this holds for all $P$, the exact diagram in
3.9 becomes
 $$
\spreadmatrixlines{3pt}
\matrix\format &\,\c\,\\
0&\to&N_2&\to&\Oh_D&\to&\Oh_\ell&\to&0\\
&&\bigcap&&\bigcap&&\Vert\\
0&\to&N_1&\to&\Oh_{D^+}&\to&\Oh_\ell&\to&0\\
&&\downarrow&&\downarrow\\
&&N_1/N_2&=&\Oh_{D^+}/\Oh_D.\\
\endmatrix
 $$
The splitting of $\Oh_D$ as a sheaf of rings follows as before. \QED\enddemo

\proclaim{{\rm4.10}\enspace Corollary} (I) Let $X$ be a nonnormal del Pezzo
surface. Then
 $$
\text{$X$ is tame }\iff \chi(\Oh_X)=1 \iff H^1(\Oh_X)=0;
 $$
and if this holds, $H^1(\Oh_X(n))=0$ for all $n$.

(II) Assume that $X$ is tame and write
$n=\deg_{\Oh_X(1)}X=\bigl(\Oh_X(1)\bigr)^2$. Then $\Oh_X(1)$ is very ample if
$n\ge3$, generated by its $H^0$ if $n=2$, and has a single transverse base
point if $n=1$.
\endproclaim

\demo{Proof of (I)} Obviously $h^0(\Oh_X)=1$ and
$h^2(\Oh_X)=h^0(\Oh_X(-1))=0$, and therefore in any case
$\chi(\Oh_X)=1-h^1(\Oh_X)$. Also
$\pi_*\Oh_Y/\Oh_X=\fie_*\Oh_C/\Oh_D$ which gives
 $$
\chi(\Oh_Y)-\chi(\Oh_X)=\chi(\Oh_C)-\chi(\Oh_D);
 $$
now each $C_i$ is a plane conic, so that $\chi(\Oh_{Y_i})=\chi(\Oh_{C_i})=1$
for each $C_i\subset Y_i$, and hence $\chi(\Oh_X)=\chi(\Oh_D)$. But the
conclusion of 4.2--3 and Proposition~4.9 is that in the tame case $D$ is
either a line, a conic, a union of linearly independent lines through a point,
or a first order neighbourhood of $\proj^1$ in $\proj^r$. These all have
$\chi(\Oh_D)=1$; for example, in the last case,
$\chi(\Oh_D)=\chi(\Oh_{\proj^1})=1$ follows from Proposition~4.9.

$N_2=(r-1)\Oh(-1)$ by Proposition~4.9, so that $\chi(N_2)=0$, and hence
$\chi(\Oh_X)=\chi(\Oh_D)=\chi(\Oh_{D\red})$. But if
$\Oh_{D\red}\subsetneq\Oh_\ell$ then $\chi(\Oh_{D\red})\le0$. Thus tame is
equivalent to $\chi(\Oh_X)=1$.

The rest of (I) is easy: using $\Oh_D=\Oh_\ell\oplus(r-1)\Oh_\ell(-1)$
together with easy arguments on the $C_i\subset Y_i$ implies that both the
maps $H^0(Y,\Oh_Y(n))\twoheadrightarrow H^0(C,\Oh_C(n)) \twoheadrightarrow
H^0(D,\fie_*\Oh_C/\Oh_D(n))$ are surjective for $n>0$, so that
$H^1(\Oh_X(n))=0$ follows from the short exact sequence
$0\to\Oh_X(n)\to\Oh_Y(n)\to\fie_*\Oh_C/\Oh_D(n)\to0$. For $n<0$, just use
duality. \enddemo

\demo{Proof of (II)} First, a similar argument to that just given shows
that $H^0(\Oh_X(1))\twoheadrightarrow H^0(\Oh_D(1))$; since $\Oh_D(1)$ is
obviously very ample, sections of $\Oh_X(1)$ embed $D$. Next, if $X$ is
reducible, $H^0(\Oh_X(1))$ embeds every component of $Y$. For either $Y_i$ is
a plane, and $H^0(\Oh_X(1))$ embeds $C_i\subset D$, hence also $Y_i$. Or $Y_i$
is not a plane, in which case $I_{C_i}\cdot\Oh_{Y_i}(1)$ is generated by its
$H^0$, which is a direct summand of $H^0(I_D\cdot\Oh_X(1))$. Finally, if $X$
is irreducible, then $\pi\:Y\to X$ composed with the rational map
corresponding to $\Oh_X(1)$ is just the linear projection from
$Y\subset\proj^{n+1}$ from a point in the plane of the conic $C$ but not on
$C$. \QED\enddemo

\subheading{{\rm4.11}\enspace The wild case} Suppose $\cha k=p$. Then in the
diagram of 3.9, $N_1/N_2\cong\Oh_\ell(-1)$, but $\Oh_{D^+}/\Oh_D$ can be
$\Oh_\ell(Np-1)$ for any $N>0$. For this, just take the derivation
$\Oh_{D^+}\to k(x)$ given by $(f,\underline g)\mapsto af'+\sum g_i$, where
$a\in k(x)$ is a rational function having poles of order exactly $n_jp$ at any
points $P_j\in\ell$ (and zeros anywhere outside the $P_j$). It's clear from
the proof of Corollary~4.10 that the corresponding surface $X$ will have
$h^1(\Oh_X)=N(p-1)$, so $\chi(\Oh_X)=1-N(p-1)$. Shepherd-Barron kindly points
out that the Picard scheme $\Pico X$ of $X$ is simply $N(p-1)$ copies of the
additive group scheme $\Bbb G_a$: in fact by deformation theory,
$H^2(\Oh_X)=0$ implies that $\Pico X$ is reduced of dimension $h^1(\Oh_X)$; and
it can't contain $\Bbb G_m$ or an Abelian variety, since then $X$ would have
cyclic etale covers of arbitrarily large order, which is absurd in view of the
concrete description of $X$ (in the irreducible case, $X$ is homeomorphic in
the etale topology to its normalisation, a rational surface).

\proclaim{{\rm4.12}\enspace Exercise} (i) Suppose that $\cha k=p\ge3,$
and that $Y$ is smooth with affine coordinates $(x,y),$ such that $C:(y^2=0)$.
In the notation of 4.4, let $h=x^{-np}\cdot h_0$ with $h_0=h_0(x)$ a unit
and $n\ge1$. Show that the local ring $\Oh_{X,P}$ is the localisation at $0$ of
 $$
k[u,v_1,\dots,v_{p-1}],\quad\text{where $u=x^p$ and
$v_i=x^{np+i}-ih_0x^{i-1}y$ for $i=1,\dots,p-1$.}
 $$
{\rm[Hint: $\Oh_X$ consists of polynomials of the form
$f(x)-hf'(x)y+g(x,y)y^2$, where $f(x)$ and $hf'(x)$ are regular, and $g(x,y)$
is arbitrary. Show that
 $$
v_iv_j-u^nv_{i+j}=\text{unit}\times x^{i+j-2}y^2,
 $$
and find a similar trick giving the monomials $x^iy^3$.]}

(ii) Take $p=3$ in (i). Show that
 $$
\spreadmatrixlines{6pt}
\matrix\format\r\ &\c\ &\l\\
\aff^2=Y&\to&X\subset\aff^3\\
\text{{\rm given by }}(x,y)&\mapsto&(x^3,x^{3n+1}-h_0y,x^{3n+2}-2h_0xy)\\
\endmatrix
 $$
is the normalisation of the hypersurface
$X:(u^{3n+2}+uv_1^3+v_2^3=0)\subset\aff^3$, a surface having cusps in
codimension $1$ along the parametrised curve $(t^3,t^{3n+1},t^{3n+2})$.

(iii) Similarly, if $\cha k=2$ then
 $$
\spreadmatrixlines{6pt}
\matrix\format\r\ &\c\ &\l\\
\aff^2=Y&\to&X\subset\aff^3\\
\text{{\rm given by }}(x,y)&\mapsto&(x^2,x^{2n+1}+y,xy^2)\\
\endmatrix
 $$
is the normalisation of $X:(u(u^{2n+1}+v^2)^2+w^2=0)\subset\aff^3$, a surface
with cusps along $(t^2,t^{2n+1},0)$. This corresponds to $h=-x^{-2n}\cdot
h_0$ in 4.4.

\endproclaim

\heading References\endheading
\parindent=0pt

[Altman--Kleiman] A. Altman and S. Kleiman, Introduction to Gro\-then\-dieck
duality, LNM {\bf 146}

[Bourbaki] N.~Bourbaki, Alg\`ebre commutative, Hermann, Paris

[del Pezzo 1] Pasquale del Pezzo, Sulle superficie dell'ordine $n$ immerse
negli
spazi di $n+1$ dimensioni, Rend\. della R. Acc\. delle Scienze Fis\. e Mat\. di
Napoli, Sept\. 1885

[del Pezzo 2] Pasquale del Pezzo, Sulle superficie dell' $n^{\roman{no}}$
ordine
immerse nello spazio di $n$ dimensioni, Rend\. del circolo matematico di
Palermo {\bf 1} (1887), 241--271

[Enriques] Federigo Enriques, Le superficie algebriche, Zanichelli,
Bo\-l\-ogna,
1949

[Grauert and Schneider] H. Grauert and M. Schneider, Komplexe Unterr\"aume und
holomorphe Vektorb\"undel vom Rang zwei, Math\. Ann\. {\bf 230} (1977), 75--90

[Grothendieck--Hakim] A. Grothendieck (written up by M. Hakim), Modules et
foncteurs dualisants, Expos\'e 4 of Cohomologie locale des faisceaux
coh\'erents et th\'eor\`emes de Lefschetz locaux et globaux, (SGA 2),
North-Holland, 1968

[Grothendieck--Hartshorne] A. Grothendieck (written up by R. Harts\-horne),
Local cohomology, LNM {\bf 41}

[Hartshorne] R. Harts\-horne, Algebraic geometry, Springer

[Horrocks] G. Horrocks, Birationally ruled surfaces without embeddings in
regular schemes, Bull\. London Math\. Soc\. {\bf 3} (1971), 57--60

[Matsumura] H. Matsumura, Commutative ring theory, C.U.P., Cambridge, 1986

[Mori] S. Mori, Three-folds whose canonical bundles are not numerically
effective, Ann\. of Math\. (2) {\bf 116} (1982), 133--176

[Mumford] D. Mumford, Appendix to Ch\.~3 of Zariski, Algebraic surfaces, 2nd
Ed., Springer, 1971

[C3-f] M. Reid, Canonical singularities, in Journ\'ees de g\'eometrie
alg\'e\-brique \linebreak
d'Angers, ed\. A. Beauville, Sijthoff and Noordhoff, Alphen 1980,
273--310

[YPG] M. Reid, Young person's guide to canonical singularities, in Algebraic
Geometry, Bowdoin 1985, Proc\. of Symposia in Pure Math\. {\bf 46}:1, 345--414,
A.M.S., 1987

[Reid] M. Reid, Local aspects of duality, expository lecture, 12 pp., approx\.
Nov\. 1987

[Semple and Roth] J. G. Semple and L. Roth, Introduction to algebraic geometry,
Oxford

[Serre] J.-P. Serre, Groupes alg\'ebriques et corps de classes, Hermann, Paris,
1959 (English translation, Springer, 1990)

[Segre] Corrado Segre, Sulle rigate razionali in uno spazio lineare
qua\-lunque,
Atti della R. Accademia delle Scienze di Torino {\bf 19} (1883--84), 355--372.
Collected works, vol\. I, pp\. 1--16

\medskip
\noindent
Address: Math Inst., Univ\. of Warwick, Coventry CV4 7AL, England

\noindent
e-mail: Miles\@Maths.Warwick.Ac.UK

\bye